\documentclass[12pt,twoside]{amsart}
\usepackage{latexsym,amsmath,amsopn,amssymb,amsthm,amsfonts}
\usepackage[T2A]{fontenc}
\usepackage[cp1251]{inputenc}
\usepackage[english]{babel}
\usepackage{graphicx}

\begin{document}
\vspace{0.5cm}
\title[Geodesics of some sub-Riemannian metrics on $SU(1,1)\times\mathbb{R}$ and $SO_0(2,1)\times\mathbb{R}$]{Geodesics and shortest arcs of some sub-Riemannian metrics on the Lie groups $SU(1,1)\times\mathbb{R}$ and $SO_0(2,1)\times\mathbb{R}$ with three-dimensional generating distributions}
\author{I.~A.~Zubareva}
\thanks{The work was carried out within the framework of the State Contract to the IM SD RAS, project FWNF--2022--0003.}
\address{Sobolev Institute of Mathematics,\newline
	4 Koptyug Av., Novosibirsk, 630090, Russia}
\email{i\_gribanova@mail.ru}

\begin{abstract}
We find geodesics, shortest arcs,  cut loci,  first conjugate loci for some left-invariant sub-Riemannian metrics on the Lie groups
$SU(1,1)\times\mathbb{R}$ and $SO_0(2,1)\times\mathbb{R}$.		
		
\end{abstract}

\maketitle

\vspace{2mm}

\section{Introduction}

\vspace{2mm}

In [1] Mubarakzyanov obtained a classification (up to isomorphism) of all four-dimensional real Lie algebras.
In [2], using some modified version of the classification,  Biggs and Remsing described all
automorphisms of Lie algebras under consideration and classified all connected Lie groups with these Lie algebras.
According to [2], the four-dimensional real Lie algebra $\mathfrak{g}_{3,6}\oplus\mathfrak{g}_1$ has a basis
$E_1$,\,$E_2$,\,$E_3$, $E_4$ satisfying the following commutation relations:
\begin{equation}
\label{abc}
[E_1,E_2]=-E_3,\quad [E_2,E_3]=E_1,\quad [E_3,E_1]=E_2,\quad [E_i,E_4]=0,\,\,i=1,2,3.
\end{equation}
Obviously, $\mathfrak{g}_{3,6}\oplus \mathfrak{g}_1\cong\mathfrak{so}(2,1)\oplus\mathbb{R}$.

It is known  (see, for example, [2]) that $\mathfrak{so}(2,1)$ is the only three-dimensional Lie algebra for which 
the corresponding universal covering group, which we denote $\tilde{A}$, is not linearizable.
Any connected Lie group with the Lie algebra $\mathfrak{so}(2,1)$ is isomorphic to either $\tilde{A}$ or an 
$n$--fold covering  $A_n$ of the restricted Lorentz group $SO_0(2,1)$, $n\geq 1$. Moreover,  $SO_0(2,1)$ is isomorphic to $PSL(2,\mathbb{R})$, 
the two-fold cover $A_2$ is isomorphic to $SL(2,\mathbb{R})$, the four-fold cover $A_4$ is isomorphic to $M_p(2,\mathbb{R})$. 

According to [2], there are six types of connected Lie groups $G$ with the Lie algebra  $\mathfrak{so}(2,1)\oplus\mathfrak{\mathbb{R}}$, namely

1) the universal covering $\tilde{A}\times\mathbb{R}$;

2)  the $n$--fold covers $(\tilde{A}\times\mathbb{R})/\exp(2n\pi\mathbb{Z}E_3)\cong A_n\times\mathbb{R}$ of the Lie group  $SO_0(2,1)\times\mathbb{R}$;

3) the Lie group $(\tilde{A}\times\mathbb{R})/\exp(\mathbb{Z}E_4)\cong\tilde{A}\times\mathbb{T}$;

4) the $n$--fold covers $(\tilde{A}\times\mathbb{R})/\exp(n\mathbb{Z}(2\pi E_3+E_4))$ of the Lie group  

$(\tilde{A}\times\mathbb{R})/\exp(\mathbb{Z}(2\pi E_3+E_4));$

5) the $n$--fold covers $(\tilde{A}\times\mathbb{R})/(\exp(2n\pi \mathbb{Z}E_3)\exp(\mathbb{Z}E_4))\cong A_n\times\mathbb{T}$
of the Lie group 

$(\tilde{A}\times\mathbb{R})/(\exp(\mathbb{Z}E_3)\exp(\mathbb{Z}E_4))\cong SO_0(2,1)\times\mathbb{T}$;

6) the $n$--fold covers  $(\tilde{A}\times\mathbb{R})/(\exp(n\mathbb{Z}(2\pi E_3+\alpha E_4))\exp(n\mathbb{Z}E_4))$
of the Lie group 

$(\tilde{A}\times\mathbb{R})/(\exp(\mathbb{Z}(2\pi E_3+\alpha E_4))\exp(\mathbb{Z}E_4))$.

\vspace{1mm}

The following proposition was proved in [3].

\vspace{1mm}

{\bf Proposition 1.} {\sl
	A three-dimensional subspace $\mathfrak{q}$ of the Lie algebra $\mathfrak{g}=\mathfrak{so}(2,1)\oplus\mathfrak{\mathbb{R}}$ generates
	$\mathfrak{g}$ if and
	only if $\mathfrak{q}\neq \mathfrak{so}(2,1)$ and the projection of $\mathfrak{q}$  to $\mathfrak{so}(2,1)$ along $\mathfrak{\mathbb{R}}$ 
	is not a two-dimensional subalgebra of  $\mathfrak{so}(2,1)$. There
	exist five nonequivalent classes of equivalent, i.e. transformable
	into each other by some automorphism of the Lie algebra $\mathfrak{g}$, these subspaces. }

\vspace{1mm}

Put $\quad\mathfrak{q}_1={\rm span}\{E_1,E_4-E_3,E_2\}$, $\mathfrak{q}_2={\rm span}\{E_1,E_4,E_2\}$, $\mathfrak{q}_3={\rm span}\{E_1,E_4,E_3\},$
$\quad\quad\quad\quad\mathfrak{q}_4={\rm span}\{E_1,E_4-E_2,E_3\},\quad 
\mathfrak{q}_5={\rm span}\{E_4+(E_2-E_3)/2,E_2+E_3,E_1\}.$

It is easy to show that each subspace  $\mathfrak{q}_i$, $i=1,\dots,5$, generates the Lie algebra $\mathfrak{so}(2,1)\oplus\mathfrak{\mathbb{R}}$,
and these subspaces are pairwise nonequivalent.

Using the results of [4], the author found geodesics of the sub-Riemannian metrics $d_i$, $i=1,2$, on
the connected Lie group  $G$ (with the unit $\rm{Id}$) with the Lie algebra $\mathfrak{so}(2,1)\oplus\mathfrak{\mathbb{R}}$
which are defined by completely
nonholonomic left-invariant distributions $D_i$ with $D_i(\rm{Id})=\mathfrak{q}_i$ and the inner product $\langle\cdot,\cdot\rangle_i$
on $D_i(\rm{Id})$ with the orthonormal basis  $E_1,E_4-E_3,E_2$ for $i=1$ and $E_1,E_4,E_2$ for $i=2$.
Applying the methods and results of [5]\,--\,[11], the author found shortest arc, the cut loci, first conjugate sets on the Lie
groups $SU(1,1)\times\mathbb{R}$ and $SO_0(2,1)\times\mathbb{R}$ equipped with these metrics. 

\section{Preliminaries}

{\bf Theorem 1.} {\sl
	Let the basis
	\begin{equation}
	\label{abc1}
	e_1=E_1,\quad e_2=E_4-E_3,\quad e_3=E_2,\quad e_4=-E_3
	\end{equation}
	be given of the Lie algebra $\mathfrak{g}=\mathfrak{so}(2,1)\oplus\mathbb{R}$, $D_1(\rm{Id})={\rm span}(e_1,e_2,e_3)$, and the inner product
	$\langle\cdot,\cdot\rangle_1$ on $D_1(\rm{Id})$ with the orthonormal basis $e_1,e_2,e_3$. 
	Then the left-invariant distribution $D_1$  on the connected Lie group $G$ with the Lie algebra $\mathfrak{g}$ is totally nonholonomic and the pair
	$(D_1(\rm{Id}),\langle\cdot,\cdot\rangle_1)$ defines the left-invariant
	sub-Riemannian metric $d_1$ on $G$. Moreover, each geodesic  $\gamma_1=\gamma_1(t)$, $t\in\mathbb{R}$, in $(G,d_1)$
	parametrized by arclength such that $\gamma_1(0)=\rm{Id}$  is a product of two one-parameter subgroups
	\begin{equation}
	\label{geod}
	\gamma_1(t)=\gamma_1(\alpha_1,\alpha_2,\alpha_3,\beta;t)=\exp(t(\alpha_1e_1+\alpha_2e_2+\alpha_3e_3+\beta e_4))\exp(-t\beta e_4),
	\end{equation}
	where $\alpha_1,\alpha_2,\alpha_3$ and $\beta$ are some  arbitrary constants such that
	\begin{equation}
	\label{norm}
	\alpha_1^2+\alpha_2^2+\alpha_3^2=1.
	\end{equation}  }

{\sl Proof.} 
It follows from (\ref{abc}) and (\ref{abc1}) that
\begin{equation}
\label{e4}
\left[e_1,e_2\right]=\left[e_1,e_4\right]=e_3,\quad [e_1,e_3]=e_4,\quad [e_2,e_3]=-[e_3,e_4]=e_1,\quad [e_2,e_4]=0.
\end{equation}
This implies the first statement of Theorem 1.	

By (\ref{e4}) and [3, Theorem~3]  every abnormal extremal of the sub-Riemannian space $(G,d_1)$  is nonstrictly
abnormal, and so it is a geodesic. According to
[4, Theorem~9], in $(G,d_1)$ every normal geodesic $\gamma_1(t)$, $t\in\mathbb{R}$, parametrized by arclength such that $\gamma_1(0)=\rm{Id}$ 
satisfies the differential equation
\begin{equation}
\label{basic}
\dot{\gamma}(t)=dl_{\gamma(t)}(u(t)),\quad u(t)=\psi_1(t)e_1+\psi_2(t)e_2+\psi_3(t)e_3,\quad |u(0)|=1,
\end{equation}
and the absolutely continuous functions $\psi_i(t)$, $i=1,\dots,4$, satisfy the system of differential equations
\begin{equation}
\label{cosystem}
\dot{\psi}_j(t)=\sum_{k=1}^{4}\sum_{i=1}^{3}C_{ij}^k\psi_i(t)\psi_k(t),\quad j=1,\dots,4.
\end{equation}
Here $C_{ij}^k$ are structure constants in the basis  $e_1,e_2,e_3,e_4$ for $\mathfrak{g}$.

Taking into account (\ref{e4}), system (\ref{cosystem}) is written as
$$\dot{\psi}_1=-\psi_3(\psi_2+\psi_4),\quad \dot{\psi}_2=0,\quad\dot{\psi}_3=\psi_1(\psi_2+\psi_4),\quad \dot{\psi}_4=0.$$
Assign an arbitrary set of initial data
$\psi_i(0)=\alpha_i$, $i=1,2,3$, and $\psi_4(0)=-\alpha_2-\beta$
2 of the system. It is easy to see that
$\psi_4(t)\equiv -\alpha_2-\beta$,
\begin{equation}
\label{psi123}
\psi_1(t)=\alpha_1\cos\beta t+\alpha_3\sin\beta t,\quad\psi_2(t)\equiv\alpha_2,\quad
\psi_3(t)=-\alpha_1\sin\beta t+\alpha_3\cos\beta t,
\end{equation}
and the condition $|u(0)|=1$ is equivalent to (\ref{norm}).

Let us prove that (\ref{geod}) is a solution of (\ref{basic}). It follows from (\ref{e4}) that
$(\rm{ad}(e_4))=e_{13}-e_{31}$, where $(f)$  denotes the matrix of a linear map $f:\mathfrak{g}\rightarrow\mathfrak{g}$ in the basis $e_1,e_2,e_3,e_4$. Here $e_{ij}$ is a $4\times 4$--matrix,  having $1$ in the $i$th row and the $j$th column while all other entries are $0$.

Differentiating (\ref{geod}) and using (\ref{basic}), (\ref{psi123}), and formula (5) in [13, Chapter~II, \S~5],  we have
$$\dot{\gamma}_1(t)=\exp(t(\alpha_1e_1+\alpha_2e_2+\alpha_3e_3+\beta e_4))(\alpha_1e_1+\alpha_2e_2+\alpha_3e_3+\beta e_4)\exp(-t\beta e_4)$$
$$+\gamma_1(t)(-\beta e_4)=\gamma_1(t)\exp(t\beta e_4)(\alpha_1e_1+\alpha_2e_2+\alpha_3e_3+\beta e_4)\exp(-t\beta e_4)+\gamma_1(t)(-\beta e_4)$$
$$=\gamma_1(t)\exp(t\beta e_4)(\alpha_1e_1+\alpha_2e_2+\alpha_3e_3)\exp(-t\beta e_4)+\gamma_1(t)(\beta e_4)+\gamma_1(t)(-\beta e_4)=$$
$$=\gamma_1(t)\cdot[\rm{Ad}(\exp(t\beta e_4))(\alpha_1e_1+\alpha_2e_2+\alpha_3e_3)]$$
$$=\gamma_1(t)\cdot[\exp(t\beta(\rm{ad}(e_4)))(\alpha_1e_1+\alpha_2e_2+\alpha_3e_3)]$$
$$=\gamma_1(t)\cdot\left[(\alpha_1\cos\beta t+\alpha_3\sin\beta t)e_1+\alpha_2 e_2+(-\alpha_1\sin\beta t+\alpha_3\cos\beta t)e_3\right]=
\gamma_1(t)u(t).$$   
\noindent $\square$

\vspace{1mm}

{\bf Theorem 2.} {\sl
	Let the basis
	\begin{equation}
	\label{abc2}
	e_1=E_1,\quad e_2=E_4,\quad e_3=E_2,\quad e_4=-E_3
	\end{equation}
	be given of the Lie algebra $\mathfrak{g}=\mathfrak{so}(2,1)\oplus\mathbb{R}$, 
	$D_2(\rm{Id})={\rm span}(e_1,e_2,e_3)$ with the inner product
	$\langle\cdot,\cdot\rangle_2$ on
	$D_2(\rm{Id})$ with the orthonormal basis
	$e_1,e_2,e_3$.  
	Then the left-invariant distribution $D_2$  on the connected Lie group $G$ with the Lie algebra $\mathfrak{g}$ is totally nonholonomic and the pair
	$(D_2(\rm{Id}),\langle\cdot,\cdot\rangle_1)$ defines the left-invariant
	sub-Riemannian metric $d_2$ on $G$. Moreover, each  geodesic  $\gamma_2=\gamma_2(t)$, $t\in\mathbb{R}$, in  $(G,d_2)$
	parametrized by arclength such that  $\gamma_2(0)=\rm{Id}$ is a product of two one-parameter subgroups
	\begin{equation}
	\label{geod2}
	\gamma_2(t)=\gamma_2(\alpha_1,\alpha_2,\alpha_3,\beta;t)=\exp(t(\alpha_1e_1+\alpha_2e_2+\alpha_3e_3+\beta e_4))\exp(-t\beta e_4),
	\end{equation}
	where $\alpha_1,\alpha_2,\alpha_3$ and $\beta$ 
	are some arbitrary constants and (\ref{norm}) holds. }

\vspace{2mm}

{\sl Proof.} 
It follows from (\ref{abc}) and (\ref{abc2}) that
\begin{equation}
\label{ee4}
\left[e_1,e_2\right]=[e_2,e_3]=[e_2,e_4]=0,\quad [e_1,e_3]=e_4,\quad [e_1,e_4]=e_3,\quad [e_3,e_4]=-e_1.
\end{equation}	
This implies the first statement of Theorem 2.

By (\ref{ee4}) и [3, Theorem 3]  every abnormal extremal of the sub-Riemannian space  $(G,d_2)$ is nonstrictly abnormal, and so it is a geodesic. Moreover, taking into account
(\ref{ee4}), system (\ref{cosystem})  is written
as
$$\dot{\psi}_1=-\psi_3\psi_4,\quad \dot{\psi}_2=0,\quad\dot{\psi}_3=\psi_1\psi_4,\quad \dot{\psi}_4=0.$$
Assign an arbitrary set of initial data $\psi_i(0)=\alpha_i$, $i=1,2,3$, and $\psi_4(0)=-\beta$ of the system. It is easy to see that $\psi_4(t)\equiv -\beta$, and the functions
$\psi_i(t)$, $i=1,2,3$,  are defined by (\ref{psi123}).
By (\ref{ee4}), $(\rm{ad}(e_4))=e_{13}-e_{31}$.
To repeat the calculations of the last paragraph of the proof of Theorem 1, we find that (\ref{geod2}) is the
solution of (\ref{basic}). This completes the proof of Theorem 2.
\noindent $\square$

\vspace{2mm}

{\bf Proposition 2.} {\sl
	For all $(\alpha_1,\alpha_2,\alpha_3)\in\mathbb{S}^2$, $\beta\in\mathbb{R}$, and $t\in\mathbb{R}$, 
	$$\gamma_1(\alpha_1,\alpha_2,\alpha_3,\beta;t)=\gamma_2(\alpha_1,\alpha_2,\alpha_3,\beta+\alpha_2;t)\exp(t\alpha_2e_4).$$ }

{\sl Proof.} By (\ref{abc1}), (\ref{geod}), (\ref{abc2}), and (\ref{geod2}),
$\gamma_2(\alpha_1,\alpha_2,\alpha_3,\beta+\alpha_2;t)\exp(t\alpha_2e_4)$
$$=\exp(t(\alpha_1E_1+\alpha_2E_4+\alpha_3E_2-(\beta+\alpha_2)E_3))
\exp(t(\beta+\alpha_2)E_3)\exp(-t\alpha_2E_3)$$
$$=\exp(t(\alpha_1E_1+\alpha_2(E_4-E_3)+\alpha_3E_2-\beta E_3))\exp(t\beta E_3)=
\gamma_1(\alpha_1,\alpha_2,\alpha_3,\beta;t).$$
\noindent $\square$

\vspace{1mm}

{\bf Proposition 3.} {\sl
	Every abnormal extremal in  $(G,d_i)$, $i=1,2$, 
	parametrized by the arclength is
	nonstrictly abnormal and is one of the two one-parameter subgroups
	$\gamma_i(0,\alpha_2,0,\beta;t)=\exp(t\alpha_2e_2)$, $\alpha_2=\pm 1$, $\beta\in\mathbb{R},$
	or their left shifts on $(G,d_i)$. }

\vspace{1mm}

{\sl Proof.}
If $\alpha_2=\pm 1$ then $\alpha_1=\alpha_3=0$ due to (\ref{norm}).  
It follows from (\ref{geod}),  (\ref{e4}), (\ref{geod2}), (\ref{ee4}),  and the equality
$\exp(X+Y)=\exp{X}\exp{Y}$ provided that $[X,Y]=0$, where $X,Y\in\mathfrak{g}$ (see, for example, [13, Chapter~II, \S~5])
that for every $\beta\in\mathbb{R}$
we have
$$\gamma_i(0,\alpha_2,0,\beta;t)=\exp(t(\alpha_2 e_2+\beta e_4))\exp(-t\beta e_4)=\exp(\alpha_2 te_2).$$
This equality and [3] imply the remaining claims of Proposition 2.
\noindent $\square$

\vspace{1mm}

{\bf Proposition 4.} {\sl 
	Let $\gamma_i(t)=\gamma_i(\alpha_1,\alpha_2,\alpha_3,\beta;t)$,  $t\in\mathbb{R}$, be a geodesic in $(G,d_i)$, $i=1,2$, defined
	by (\ref{geod}) and (\ref{geod2}). Then for every $t_0\in\mathbb{R}$ we have
	$$\gamma_i(t_0)^{-1}\gamma_i(t)=\gamma_i(\alpha_1\cos\beta t_0+\alpha_3\sin\beta t_0,\alpha_2,-\alpha_1\sin\beta t_0+\alpha_3\cos\beta t_0,\beta;t-t_0).$$ }

{\sl Proof.} 
Using (\ref{geod}), (\ref{geod2}), and the equality $({\rm ad}(e_4))=e_{13}-e_{31}$ in the proofs of Theorems 1, 2, and (5) в [9, Chapter~II, \S~5], we get
$$\gamma_i(t_0)^{-1}\gamma_i(t)=\exp(t_0\beta e_4)\exp(-t_0(\alpha_1e_1+\alpha_2e_2+\alpha_3e_3+\beta e_4))$$
$$\times\exp(t(\alpha_1e_1+\alpha_2e_2+\alpha_3e_3+\beta e_4))\exp(-t\beta e_4)$$
$$=\exp(t_0\beta e_4)\exp((t-t_0)(\alpha_1e_1+\alpha_2e_2+\alpha_3e_3+\beta e_4))\exp(-t_0\beta e_4)\exp(-(t-t_0)\beta e_4)$$
$$=\exp[\rm{Ad}(\exp(t_0\beta e_4))((t-t_0)(\alpha_1e_1+\alpha_2e_2+\alpha_3e_3+\beta e_4))]\cdot\exp(-(t-t_0)\beta e_4)$$
$$=\exp[\exp(\rm{ad}(t_0\beta e_4))((t-t_0)(\alpha_1e_1+\alpha_2e_2+\alpha_3e_3+\beta e_4))]\cdot\exp(-(t-t_0)\beta e_4)$$
$$=\exp((t-t_0)((\alpha_1\cos\beta t_0+\alpha_3\sin\beta t_0)e_1+\alpha_2e_2$$
$$+(-\alpha_1\sin\beta t_0+\alpha_3\cos\beta t_0)e_3+\beta e_4)\cdot\exp(-(t-t_0)\beta e_4)$$
$$=\gamma_i(\alpha_1\cos\beta t_0+\alpha_3\sin\beta t_0,\alpha_2,-\alpha_1\sin\beta t_0+\alpha_3\cos\beta t_0,\beta;t-t_0).$$
\noindent $\square$

\section{Sub-Riemannnian geodesics on $SU(2)\times\mathbb{R}$}

Recall that the Lie group $SU(1,1)$  consists of all complex $2\times 2$--matrices 
with determinant $1$, preserving the indefinite Hermitian form $|z_1|^2-|z_2|^2$ in $\mathbb{C}^2$:
$$SU(1,1)=\left\{(A,B):=\left(\begin{array}{cc}
A & B \\
\overline{B} & \overline{A}
\end{array}\right)\left|\,\,A,\,B\in\mathbb{C},\,\,|A|^2-|B|^2=1\right.\right\}.$$

The Lie algebra $\mathfrak{su}(1,1)$ of the Lie group $SU(1,1)$ is given by
$$\mathfrak{su}(1,1)=\left\{\left(\begin{array}{cc}
iX & Y \\
\overline{Y} & -iX
\end{array}\right)\left|\,\,X\in\mathbb{R},\,Y\in\mathbb{C}\right.\right\},\quad \mathfrak{su}(1,1)\cong\mathfrak{so}(2,1).$$

The trivial abelian extension
\begin{equation}
\label{defSU(1,1)xR}
SU(1,1)\times\mathbb{R}=\left\{(A,B,v):=\left(\begin{array}{ccc}
A & B & 0 \\
\overline{B} & \overline{A} & 0 \\
0 & 0 & e^v
\end{array}\right)\left|\,A,\,B\in\mathbb{C},\,\,|A|^2-|B|^2=1;\,\,v\in\mathbb{R}\right.\right\}
\end{equation}
of $SU(1,1)$ is the Lie group whose the Lie algebra $\mathfrak{su}(1,1)\oplus\mathfrak{\mathbb{R}}$ has some basis $E_1,E_2,E_3,E_4$ satisfying (\ref{abc}):
\begin{equation}
\label{Basis1}
E_1=\frac{1}{2}(e_{12}+e_{21}),\quad E_2=\frac{i}{2}(e_{12}-e_{21}),\quad E_3=\frac{i}{2}(e_{11}-e_{22}),\quad E_4=e_{33}.
\end{equation}
Here $e_{ij}$, $i,j=1,\dots,3$, denotes the $3\times 3$--matrix, 
having $1$ in the $i$th row and the $j$th column while all
other entries are $0$.

\vspace{2mm}

{\bf Remark 1.} {\sl 
	The Lie group $SU(1,1)$ is isomorphic to the special linear Lie group $SL(2,\mathbb{R})$
	consisting of all real $2\times 2$--matrices with determinant $1$.
	The desired isomorphism $\Omega:SU(1,1)\rightarrow SL(2,\mathbb{R})$ is given by formula
	$$\Omega(g)=hgh^{-1},\,\, h=\left(\begin{array}{cc}
	1 & i \\
	i & 1
	\end{array}\right);\quad \Omega(A,B)=\left(\begin{array}{cc}
	{\rm Re}(A)+{\rm Im}(B) & {\rm Im}(A)+{\rm Re}(B) \\
	{\rm Re}(B)-{\rm Im}(A) & {\rm Re}(A)-{\rm Im}(B)
	\end{array}\right),$$
	can naturally be extended to the isomorphism
	$\tilde{\Omega}:SU(1,1)\times\mathbb{R}\rightarrow SL(2,\mathbb{R})\times\mathbb{R}$.
	
	It is easy to check that the differential $d\tilde{\Omega}({\rm Id})$
	is an isomorphism of Lie algebras
	$\mathfrak{su}(1,1)\oplus\mathfrak{\mathbb{R}}$ and
	$\mathfrak{sl}(2)\oplus\mathfrak{\mathbb{R}}$
	translating the basis $E_1,E_2,E_3,E_4$ of Lie algebra $\mathfrak{su}(1,1)\oplus\mathfrak{\mathbb{R}}$ given by (\ref{Basis1}),
	to the same name basis of Lie algebra $\mathfrak{sl}(2)\oplus\mathfrak{\mathbb{R}}$
	given by
	$$E_1=\frac{1}{2}(e_{12}+e_{21}),\quad E_2=\frac{1}{2}(e_{11}-e_{22}),\quad  E_3=\frac{1}{2}(e_{12}-e_{21}),\quad E_4=e_{33}.$$ }

\vspace{1mm}

{\bf Theorem 3.} {\sl
	Let
	\begin{equation}
	\label{mn01}
	m_2=t/2,\quad n_2=1,\quad\text{if  }\,\,\beta^2=1-\alpha_2^2;
	\end{equation}
	\begin{equation}
	\label{mn02}
	m_2=\frac{\sinh\frac{t\sqrt{1-\alpha_2^2-\beta^2}}{2}}{\sqrt{1-\alpha_2^2-\beta^2}},\quad
	n_2=\cosh\frac{t\sqrt{1-\alpha_2^2-\beta^2}}{2},\quad\text{if  }\,\,
	\beta^2<1-\alpha_2^2;
	\end{equation}
	\begin{equation}
	\label{mn03}
	m_2=\frac{\sin\frac{t\sqrt{\beta^2+\alpha_2^2-1}}{2}}{\sqrt{\beta^2+\alpha_2^2-1}},\quad n_2=\cos\frac{t\sqrt{\beta^2+\alpha_2^2-1}}{2},\quad\text{if  }\,\,
	\beta^2>1-\alpha_2^2.
	\end{equation}
	
	If $\alpha_2\neq\pm 1$ then the geodesic $\tilde{\gamma}_2(t)=\tilde{\gamma}_2(\alpha_1,\alpha_2,\alpha_3,\beta;t)=(A,B,v)(t)$, $t\in\mathbb{R}$,  of the sub-Riemannian space $(SU(1,1)\times\mathbb{R},d_2)$
	defined by the equality (\ref{geod2}) is given by the formulas
	\begin{equation}
	\label{form1}
	A=\left(n_2\cos\frac{\beta t}{2}+\beta m_2\sin\frac{\beta t}{2}\right)+
	\left(n_2\sin\frac{\beta t}{2}-\beta m_2\cos\frac{\beta t}{2}\right)i,
	\end{equation}
	\begin{equation}
	\label{fform2}
	B=m_2\sqrt{1-\alpha_2^2}\left[\cos\left(\frac{\beta t}{2}+\varphi_0\right)-\sin\left(\frac{\beta t}{2}+\varphi_0\right)i\right],\quad v=\alpha_2t,
	\end{equation}
	where 
	\begin{equation}
	\label{varphi}
	\cos\varphi_0=\frac{\alpha_1}{\sqrt{1-\alpha_2^2}},\quad \sin\varphi_0=-\frac{\alpha_3}{\sqrt{1-\alpha_2^2}}.
	\end{equation}
	
	If $\alpha_2=\pm 1$ then the geodesic of the sub-Riemannian space
	$(SU(1,1)\times\mathbb{R},d_2)$, defined by the formula (\ref{geod2}), is equal to
	$\tilde{\gamma}_2(t)=(e,\alpha_2t)$, $t\in\mathbb{R}$, where $e$ is the identity $2\times 2$--matrix. }

\vspace{1mm}

{\sl Proof.}
Let $\alpha_2\neq \pm 1$. Using (\ref{norm}), (\ref{abc2}), (\ref{Basis1}), and (\ref{mn01})--(\ref{mn03}), 
it is easy to get that
$$\exp(t(\alpha_1e_1+\alpha_2e_2+\alpha_3e_3+\beta e_4))=
\left(\begin{array}{ccc}
n_2-\beta m_2i & (\alpha_1+\alpha_3i)m_2 & 0 \\
(\alpha_1-\alpha_3i)m_2 & n_2+\beta m_2 i & 0 \\
0 & 0 & e^{\alpha_2t}
\end{array}\right),$$
$$\exp(-t\beta e_4)=
\left(\begin{array}{ccc}
\cos\frac{\beta t}{2}+i\sin\frac{\beta t}{2} & 0 & 0  \\
0 & \cos\frac{\beta t}{2}-i\sin\frac{\beta t}{2} & 0 \\
0 & 0 & 1
\end{array}\right).$$
By (\ref{geod2}), it remains to multiply these matrices.

The statement  of Theorem 3 for $\alpha_2=\pm 1$ follows from Propositions 3, (\ref{abc2}), (\ref{Basis1}).
\noindent $\square$

\vspace{2mm}

The following is immediate from Theorem 3 and Propositions 2, 3.

\vspace{2mm}

{\bf Theorem 4.} {\sl
	Let
	$$m_1=\frac{t}{2},\quad n_1=1,\quad\text{if  }\,\, (\alpha_2+\beta)^2=1-\alpha_2^2;$$
	$$m_1=\frac{\sinh\frac{t\sqrt{1-\alpha_2^2-(\alpha_2+\beta)^2}}{2}}{\sqrt{1-\alpha_2^2-(\alpha_2+\beta)^2}},\quad n_1=\cosh\frac{t\sqrt{1-\alpha_2^2-(\alpha_2+\beta)^2}}{2},$$
	if $(\alpha_2+\beta)^2<1-\alpha_2^2$;
	$$m_1=\frac{\sin\frac{t\sqrt{(\alpha_2+\beta)^2+\alpha_2^2-1}}{2}}{\sqrt{(\alpha_2+\beta)^2+\alpha_2^2-1}},\quad n_1=\cos\frac{t\sqrt{(\alpha_2+\beta)^2+\alpha_2^2-1}}{2},$$
	if $(\alpha_2+\beta)^2>1-\alpha_2^2$.
	
	If $\alpha_2\neq\pm 1$  then the geodesic $\tilde{\gamma}_1(t)=\tilde{\gamma}_1(\alpha_1,\alpha_2,\alpha_3,\beta;t)=(A,B,v)(t)$, with $t\in\mathbb{R}$,  of the sub-Riemannian space $(SU(1,1)\times\mathbb{R},d_1)$ defined by the equality (\ref{geod}) is given by the formulas
	$$A=\left(n_1\cos\frac{\beta t}{2}+(\alpha_2+\beta)m_1\sin\frac{\beta t}{2}\right)+
	\left(n_1\sin\frac{\beta t}{2}-(\alpha_2+\beta)m_1\cos\frac{\beta t}{2}\right)i,$$
	$$B=m_1\sqrt{1-\alpha_2^2}\left[\cos\left(\frac{\beta t}{2}+\varphi_0\right)-\sin\left(\frac{\beta t}{2}+\varphi_0\right)i\right],\quad v=\alpha_2t,$$
	and (\ref{varphi}) holds.
	
	If $\alpha_2=\pm 1$ then the geodesic $\tilde{\gamma}_1(t)=(A,B,v)(t)$, with $t\in\mathbb{R}$, of the sub-Riemannian space $(SU(1,1)\times\mathbb{R},d_1)$ defined by the equality (\ref{geod}) is given by the formulas
	\begin{equation}
	\label{form3}
	A=e^{-i(\alpha_2t/2)},\quad B\equiv 0,\quad v=\alpha_2t.
	\end{equation}  }

{\bf Remark 2.} {\sl
	In the case $\alpha_2\neq\pm 1$ we will also denote the geodesic $\tilde{\gamma}_i(\alpha_1,\alpha_2,\alpha_3,\beta;t)$ by
	$\tilde{\gamma}_i(\varphi_0,\alpha_2,\beta;t)$, if (\ref{varphi}) holds. }

\section{The Cut Locus and First Conjugate Locus in $SU(1,1)\times\mathbb{R}$}

In this section, we will recall necessary definitions from the paper [10] by Sachkov for the Lie group $G$ considered in Section 2.

The time moment $\widehat{t}>0$  is called the conjugate time for the normal geodesic
$\gamma_i(t)={\rm Exp}_i(\lambda,t)$, $\lambda=(\alpha_1,\alpha_2,\alpha_3,\beta)\in C:=\mathbb{S}^2\times\mathbb{R}$, of the sub-Riemannian space
$(G,d_i)$ given by the formula (\ref{geod}) for $i=1$ or (\ref{geod2}) for $i=2$,
if $(\lambda,\widehat{t})$  is the critical point of the exponential map, i.e., the differential
$$\left({\rm Exp}_i\right)_{\ast(\lambda,\widehat{t})}:T_{(\lambda,\widehat{t})}(C\times\mathbb{R}_{+})\rightarrow T_{\widehat{g}_i}G,
\quad\text{where  }\widehat{g}_i={\rm Exp}_i(\lambda,\widehat{t}),$$
is degenerate. The first conjugate time along the normal geodesic $\gamma_i(t)$ is defined by
$$t^1_{{\rm conj}}=\inf\{t>0\mid\,t \text{ is a conjugate time along } \gamma_i(\cdot)\}.$$
If $t^1_{{\rm conj}}>0$, then $\gamma_i(t^1_{{\rm conj}})$  is called the first conjugate point.

The first conjugate locus of the sub-Riemannian space $(G,d_i)$  is the set ${\rm Conj}^{1}_i$ 
of all first conjugate points along normal geodesics $\gamma_i(t)$ starting at  ${\rm Id}$.  
The cut locus  (for ${\rm Id}$) of the sub-Riemannian space $(G,d_i)$  is the set ${\rm Cut}_i$  of the endpoints $g\in G$ of all shortest arcs starting at ${\rm Id}$ and noncontinuable beyond $g$.

It's easy to show, given Proposition 3, that $t^1_{{\rm conj}}=0$ along a nonstrictly abnormal extremal $\gamma_i(t)={\rm exp}(\pm t e_2)$, $t\in\mathbb{R}$.

\vspace{1mm}

{\bf Proposition 5.} {\sl
	Let
	\begin{equation}
	\label{w12}
	w_1=\sqrt{(\alpha_2+\beta)^2+\alpha_2^2-1},\quad w_2=\sqrt{\beta^2+\alpha_2^2-1},
	\quad w_1,w_2\in\mathbb{R}_{+}\cup \{0\}\cup \left(i\mathbb{R}_{+}\right).
	\end{equation}
	The time moment $\widehat{t}>0$ is the conjugate time along the geodesic
	$\tilde{\gamma}_i(t)=\tilde{\gamma}_i(\varphi_0,\alpha_2,\beta;t)$,  $t\in\mathbb{R}$,  
	of the sub-Riemannian space
	$(SU(1,1)\times\mathbb{R},d_i)$, $i=1,2$, if and only if $w_i>0$ and
	$$\sin\frac{w_i\widehat{t}}{2}\left(\sin\frac{w_i\widehat{t}}{2}-\frac{w_i\widehat{t}}{2}\cos\frac{w_i\widehat{t}}{2}\right)=0.$$  }

{\sl Proof.}
Let $i=2$. For convenience of notation, we  omit the subscript of the functions
$\tilde{\gamma}_2(t)$,\,$w_2$,
$n_2(t)$, and $m_2(t)$.
By (\ref{defSU(1,1)xR}) and (\ref{Basis1}), 
$$\tilde{\gamma}(t)=2B_1(t)E_1+2B_2(t)E_2+2A_2(t)E_3+e^{\alpha_2t}E_4+2A_1(t)E_5,$$
where $E_5=\frac{1}{2}(E-E_4)$ and  $E$ is the identity $3\times 3$--matrix, and
\begin{equation}
\label{AB}
A_1={\rm Re}(A),\quad A_2={\rm Im}(A),\quad B_1={\rm Re}(B),\quad B_2={\rm Im}(B).
\end{equation}

Let $w\neq 0$. Arguing in the same way as in the proof of Proposition 4 in [12], we obtain that if the vectors
$\tilde{\gamma}^{\prime}_{t}$,
$\tilde{\gamma}^{\prime}_{\alpha_2}$, $\tilde{\gamma}^{\prime}_{\beta}$,
$\tilde{\gamma}^{\prime}_{\varphi_0}$ are linearly dependent for $t\neq 0$, then the equality
\begin{equation}
\label{AB2}
\frac{m(2m-nt)(1-\alpha_2^2)t}{4w^2}=0
\end{equation}
holds. Hence it follows that  $m=0$ or $m=nt/2$. If $m=0$ then
$\tilde{\gamma}^{\prime}_{\varphi_0}=0$; if $m=nt/2$ then
$\tilde{\gamma}'_{\beta}\parallel\tilde{\gamma}'_{\varphi_0}.$
Thus, in both cases, $\tilde{\gamma}^{\prime}_{t}$,
$\tilde{\gamma}^{\prime}_{\alpha_2}$, $\tilde{\gamma}^{\prime}_{\beta}$,
$\tilde{\gamma}^{\prime}_{\varphi_0}$ are linearly dependent.

Note that for $w^2<0$ the equality (\ref{AB2}) holds only for $t=0$, i.e.,
the geodesic $\tilde{\gamma}(t)$, $t\in\mathbb{R}$, does not contain conjugate points for $w^2<0$.

Let now $w=0$. It is easy to see that for $t\neq 0$ the vectors $\tilde{\gamma}^{\prime}_{t}$,
$\tilde{\gamma}^{\prime}_{\alpha_2}$, $\tilde{\gamma}^{\prime}_{\varphi_0}$ 
are linearly dependent if and only if the vectors $t\tilde{\gamma}^{\prime}_{t}-\alpha_2\tilde{\gamma}^{\prime}_{\alpha_2}$ and $\tilde{\gamma}^{\prime}_{\varphi_0}$ 
are collinear. In this case, it is necessary that
\begin{equation}
\label{h}
t\left[(B_1)'_{\varphi_0}(B_2)'_{t}-(B_1)'_{t}(B_2)'_{\varphi_0}\right]-\alpha_2\left[(B_1)'_{\varphi_0}(B_2)'_{\alpha_2}-
(B_1)'_{\alpha_2}(B_2)'_{\varphi_0}\right]=0.
\end{equation}
It follows from  (\ref{mn01}) and (\ref{fform2}) that
$(B_1)'_{\varphi_0}=B_2$ and $(B_2)'_{\varphi_0}=-B_1$. 
Therefore, the equality (\ref{h}) can be written as
$$t(B_1^2+B_2^2)'_{t}-\alpha_2(B_1^2+B_2^2)'_{\alpha_2}=0.$$
Due to (\ref{fform2}), $B_1^2+B_2^2=t^2(1-\alpha_2^2)/4$. We get
$$0=t\left(t^2(1-\alpha_2^2)/4\right)'_{t}-\alpha_2\left(t^2(1-\alpha_2^2)/4\right)'_{\alpha_2}=t^2/2.$$
Consequently,  the geodesic $\tilde{\gamma}(t)$, $t\in\mathbb{R}$, does not contain conjugate points for $w=0$.

The case $i=1$ is considered similarly.
\noindent $\square$

\vspace{2mm}

{\bf Corollary 1.} {\sl 
	1. The $n$--conjugate time $t^{n}_{\rm{conj}}$ along the geodesic
	$\tilde{\gamma}_i(t)=\tilde{\gamma}_i(\varphi_0,\alpha_2,\beta;t)$, where $w_i>0$, 
	of the sub-Riemannian space $(SU(1,1)\times\mathbb{R},d_i)$, $i=1,2$, has the form
	$$t^{2m-1}_{\rm{conj}}=\frac{2\pi m}{w_i},\quad t^{2m}_{\rm{conj}}=\frac{2x_m}{w_i},\quad m\in \mathbb{N},$$
	where $\{x_1,x_2,\dots\}$ are the ascending positive roots of the equation
	$\tan x=x$.
	
	2. The first conjugate locus ${\rm Conj}_i^{1}$ of the sub-Riemannian space $(SU(1,1)\times\mathbb{R},d_i)$, $i=1,2$, has the form
	$${\rm Conj}_i^{1}=\{(A,B,v)\in SU(1,1)\times\mathbb{R}\left|\,|A|=1,\,B=0,\,v\in\mathbb{R}\right.\}.$$ }

{\sl Proof.} The first statement of Corollary 1 follows from Proposition 5. 
Further, by Theorems 3 and 4,
\begin{equation}
\label{conj}
\tilde{\gamma}_j\left(t^1_{{\rm conj}}\right)=\tilde{\gamma}_j\left(\frac{2\pi}{w_j}\right)=\left(-e^{i(\beta\pi/w_{j})},0,\frac{2\alpha_2\pi}{w_j}\right),\quad j=1,2.
\end{equation}

Consider the functions
$$f_j(\alpha_2,\beta)=\frac{\beta\pi}{w_{j}},\quad v_j(\alpha_2,\beta)=\frac{2\alpha_2\pi}{w_{j}},\quad 
(\alpha_2,\beta)\in (-1,1)\times\mathbb{R},\,\,w_j>0,\quad j=1,2.$$ 

It follows from (\ref{w12})  that ranges of $f_2(\alpha_2,\beta)$ и $v_2(\alpha_2,\beta)$
are $\mathbb{R}\backslash [-\pi,\pi]$ and $\mathbb{R}$ respectively, and
$$v_2(\alpha_2,\beta)=\frac{2\alpha_2\sqrt{f_2^2(\alpha_2,\beta)-\pi^2}}{\sqrt{1-\alpha_2^2}}.$$ 
From here and (\ref{conj}) we obtain ${\rm Conj}_2^{1}$.

By (\ref{w12}) the range of each of the functions
$f_1(\alpha_2,\beta)$, $v_1(\alpha_2,\beta)$ is  $\mathbb{R}$, and
$$f_1(\alpha_2,\beta)+\frac{1}{2}v_1(\alpha_2,\beta)=f_2(\alpha_2,\alpha_2+\beta)\in\mathbb{R}\backslash [-\pi,\pi].$$
From here we obtain ${\rm Conj}_1^{1}$.
\noindent $\square$

\vspace{2mm}

{\bf Proposition 6.} {\sl 
	\label{ss1}
	If $\alpha_2=\pm 1$ 1 then the geodesic $\tilde{\gamma}_i(t)=\exp(t\alpha_2 e_2)$, $t\in\mathbb{R}$, 
	is a metric line of the sub-Riemannian space $(SU(1,1)\times\mathbb{R},d_i)$, $i=1,2$; i.e., the segment 
	$\tilde{\gamma}_i(t)$, $t_0\leq t\leq t_1$, is a shortest arc for all  $t_0,t_1\in\mathbb{R} $. }

\vspace{2mm}

{\sl Proof.}
It follows from Theorems 3 and 4 that  the length of the arc of the abnormal extremal
in $(SU(1,1)\times\mathbb{R},d_i)$, $i=1,2$, connecting
${\rm Id}$ and $(A,B,v)$, $0\neq v\in\mathbb{R}$, is equal to $|v|$, and 
the length of the arc of the normal geodesic
$\tilde{\gamma}_i(\varphi_0,\alpha_2,\beta;t)$ connecting these elements, is equal to $|v/\alpha_2|>|v|$.
\noindent $\square$

\vspace{2mm}

{\bf Proposition 7.} {\sl 
	1. The geodesic
	$\tilde{\gamma}_1(t)=\tilde{\gamma}_1(\varphi_0,\alpha_2,-\alpha_2;t)$, $t\in\mathbb{R}$, in the
	sub-Riemannian space $(SU(1,1)\times\mathbb{R},d_1)$ is a metric line.	
	
	2. The geodesic
	$\tilde{\gamma}_2(t)=\tilde{\gamma}_2(\varphi_0,\alpha_2,0;t)$, $t\in\mathbb{R}$, in the sub-Riemannian space
	$(SU(1,1)\times\mathbb{R},d_2)$ is a metric line.}
\vspace{1mm}

{\sl Proof.}
Let us prove item~2 of Proposition 7.
Let $\alpha_2^{\ast}\neq\pm 1$, $\varphi_0^{\ast}\in\mathbb{R}$, $T>0$. Put $(A,B,v):=\tilde{\gamma}_2(\varphi_0^{\ast},\alpha_2^{\ast},0;T)$.
It follows from Theorem~3 that
\begin{equation}
\label{cc1}
\left\{
\begin{array}{l}
A=\cosh\frac{T\sqrt{1-\alpha_2^{\ast 2}}}{2}, \\
|B|=\sinh\frac{T\sqrt{1-\alpha_2^{\ast 2}}}{2}, \\
 v=\alpha_2^{\ast}T
\end{array}\right.
\quad\Rightarrow\quad
T=\sqrt{v^2+4\left(\sinh^{-1}|B|\right)^2}.
\end{equation}

Assume that ${\rm Id}$ and $(A,B,v)$ can be connected by another geodesic of the sub-Riemannian space
$(SU(1,1)\times\mathbb{R},d_2)$, i.e.,
$(A,B,v)=\tilde{\gamma}_2(\varphi_0,\alpha_2,\beta;t_0)$. Then $\beta\neq 0$ by (\ref{cc1}) and (\ref{fform2}).  By Theorem~3,
\begin{equation}
\label{cc2}
|B|=|m_2(t_0)|\sqrt{1-\alpha_2^2},\quad v=\alpha_2t_0.
\end{equation}
Due to (\ref{cc1}) and (\ref{cc2}),
\begin{equation}
\label{cc3}
|t_0|>T\quad\Leftrightarrow\quad |t_0|\sqrt{1-\alpha_2^2}>2\sinh^{-1}|B|.
\end{equation}
\vspace{1mm}

{\bf Lemma 1.} {\sl 
	Let $|B|>0$, $0<\mu<1$, and $0<\eta\leq 1/|B|$. Then 
	$$\frac{|B|}{\sqrt{1+|B|^2}}<\sinh^{-1}|B|<|B|,$$
	$$\quad \sinh^{-1}(\mu |B|)>\mu\cdot\sinh^{-1}|B|,\quad 
	\sin^{-1}(\eta|B|)>\eta\cdot\sinh^{-1}|B|.$$ }

{\sl Proof.}
To prove the first inequality, it suffices to note that the functions
$$f_1(\varkappa)=\varkappa-\sinh^{-1}\varkappa,\quad 
f_2(\varkappa)=\sinh^{-1}\varkappa-\frac{\varkappa}{\sqrt{1+\varkappa^2}}$$
increase on the ray $[0,+\infty)$ and $f_1(0)=f_2(0)=0$. 

To prove the second inequality, it suffices to note that the function
$$F_2(\varkappa)=\sinh^{-1}(\varkappa|B|)-\varkappa \sinh^{-1}|B|,\quad\varkappa\in [0,1],$$
increases on $[0,\varkappa_0]$ and decreases on $[\varkappa_0,1]$, where
$$\varkappa_0=\sqrt{\frac{1}{(\sinh^{-1}|B|)^2}-\frac{1}{|B|^2}}\in (0,1),\quad F_2(0)=F_2(1)=0.$$

To prove the third inequality, it suffices to note that the function
$$F_3(\varkappa)=\sin^{-1}(\varkappa|B|)-\varkappa\sinh^{-1}|B|,\quad\varkappa\in [0,1/|B|),$$
is increasing, and $F_3(0)=0$, because 
$$F_3^{\prime}(\varkappa)=\frac{|B|}{\sqrt{1-\varkappa^2|B|^2}}-\sinh^{-1}|B|>|B|-\sinh^{-1}|B|>0.$$
Lemma~1 is proved.
\noindent $\square$

\vspace{1mm}

Let $\beta^2=1-\alpha_2^2$. By  (\ref{mn01}) and (\ref{cc2}), $|t_0|=\frac{2|B|}{\sqrt{1-\alpha_2^2}}$. It follows from here, (\ref{cc3}) and Lemma 1 that $|t_0|>T$.

Let $\beta^2<1-\alpha_2^2$.  Let us denote  $\varkappa=\sqrt{1-\alpha_2^2-\beta^2}$. It follows from (\ref{mn02}), (\ref{cc2}), and (\ref{cc3}) that
$$|t_0|=\frac{2}{\varkappa}\sinh^{-1}\left(\frac{\varkappa |B|}{\sqrt{\varkappa^2+\beta^2}}\right)>T\quad\Leftrightarrow\quad
\sinh^{-1}\left(\frac{\varkappa |B|}{\sqrt{\varkappa^2+|B|^2}}\right)>\frac{\varkappa\sinh^{-1}|B|}{\sqrt{\varkappa^2+|B|^2}}.$$
The last inequality is satisfied by virtue of Lemma 1.

Let $\beta^2>1-\alpha_2^2$.  Let us denote $\varkappa=\sqrt{\beta^2+\alpha_2^2-1}$.  It follows from (\ref{mn03}), (\ref{cc2}), and (\ref{cc3}) that
$$|t_0|=\frac{2}{\varkappa}\sin^{-1}\left(\frac{\varkappa |B|}{\sqrt{\beta^2-\varkappa^2}}\right)>T\quad\Leftrightarrow\quad
\sin^{-1}\left(\frac{\varkappa |B|}{\sqrt{\beta^2-\varkappa^2}}\right)>\frac{\varkappa\sinh^{-1}|B|}{\sqrt{\beta^2-\varkappa^2}}.$$
The last inequality is satisfied by virtue of Lemma 1.

The first statement of Proposition 7 is proved in a similar way.
\noindent $\square$

\vspace{2mm}

We need the following proposition from [11].

\vspace{2mm}

{\bf Proposition 8.} {\sl 
	If in the Lie group with left-invariant sub-Riemannian metric two points are joined
	by two different normal geodesics parametrized by arclength of equal length, then neither of these
	geodesics is a shortest arc nor a part of a longer shortest arc.}

\vspace{2mm}

{\bf Proposition 9.} {\sl 
	If the geodesic segment
	$\tilde{\gamma}_i(t)=\tilde{\gamma}_i(\varphi_0,\alpha_2,\beta;t)$, $0\leq t\leq T$, of the sub-Riemannian space $(SU(1,1),d_i)$, $i=1,2$, is a shortest arc and $w_i>0$, then $T\leq 2\pi/w_i$. }

\vspace{2mm}

{\sl Proof.} It follows from (\ref{conj}) that $\tilde{\gamma}_i(\varphi_0,\alpha_2,\beta;2\pi/w_i)$ 
does not depend on $\varphi_0$. 
It remains to apply Proposition 8.
\noindent $\square$

\vspace{1mm}

{\bf Proposition 10.} {\sl 
	The cut locus ${\rm Cut}_2$ of the sub-Riemannian space 
	$(SU(1,1)\times\mathbb{R},d_2)$,  corresponding to ${\rm Id}$, is the set
	$${\rm Cut}_2={\rm Cut}_2^{{\rm loc}}\cup{\rm Cut}_2^{{\rm glob}},$$
	where 
	$${\rm Cut}_2^{{\rm loc}}=\{(A,B,v)\in SU(1,1)\times\mathbb{R}\mid  A\neq 1,\,\,B=0,\,\,v\in\mathbb{R}\},$$
	$${\rm Cut}_2^{{\rm glob}}=\{(A,B,v)\in SU(1,1)\times\mathbb{R}\mid {\rm Re}(A)<-1,\,\,{\rm Im}(A)=0,\,\,v\in\mathbb{R}\}.$$ }

{\sl Proof.} 
Let us prove that ${\rm Cut}_2^{{\rm loc}}\subset{\rm Cut}_2$. 
Let $g=(A,B,v)\in {\rm Cut}_2^{{\rm loc}}$. By Theorem 3, Proposition 6 and the Cohn–Vossen Theorem about the existence of shortest arcs there exists the  shortest arc 
$\tilde{\gamma}_2(\varphi_0,\alpha_2,\beta;t)$, $0\leq t\leq T$, joining ${\rm Id}$ and $g$. Since $B=0$, then $m_2(T)=0$ by Theorem 3. Then 
$\beta^2>1-\alpha_2^2$ and $T=\frac{2\pi}{w_2}$ (see Proposition 9); moreover, by virtue of (\ref{conj}),
$\tilde{\gamma}_2(\varphi_0,\alpha_2,\beta;T)=\tilde{\gamma}_2(\varphi,\alpha_2,\beta;T)$
for any $\mathbb{\varphi}\in\mathbb{R}$.
Then $g\in {\rm Cut}_2$ by virtue of Proposition 8.

Let us prove that ${\rm Cut}_2^{{\rm glob}}\subset{\rm Cut}_2$. Let $g=(A,B,v)\in {\rm Cut}_2^{{\rm glob}}$ and
$\tilde{\gamma}_2(\varphi_0,\alpha_2,\beta;t)$, $0\leq t\leq T$,
is the shortest arc joining  ${\rm Id}$ and $g$. It follows from (\ref{form1}) and (\ref{cc1}) that
$A=-\sqrt{n_2^2(T)+\beta^2m_2^2(T)}$, $\beta\neq 0$. It follows from here, (\ref{fform2}) and the oddness of the function $m_2(t)$ that
$$\tilde{\gamma}_2(\varphi_0,\alpha_2,\beta;T)=\tilde{\gamma}_2(\beta T+\varphi_0+\pi,\alpha_2,\beta;-T).$$
Therefore $g\in {\rm Cut}_2$ by Proposition 8.

It follows from Theorem 3 and Propositions 6, 8 that $(1,0,v)\notin {\rm Cut}_2$, $v\in\mathbb{R}$.
We assume that in $(SU(1,1),d_2)$ for some $g=(A,B,v)\notin {\rm Cut}_2^{{\rm loc}}\cup{\rm Cut}_2^{{\rm glob}}$
there exist two distinct shortest arcs
$\tilde{\gamma}_2(\varphi_0,\alpha_2,\beta;t)$, $\tilde{\gamma}_2(\varphi_0^{\ast},\alpha_2^{\ast},\beta^{\ast};t)$, $0\leq t\leq T$, joining ${\rm Id}$ and $g$. By Proposition 7, 
$\beta\neq 0$ and $\beta^{\ast}\neq 0$. Since $v=\alpha_2T=\alpha_2^{\ast}T$, then $\alpha_2=\alpha_2^{\ast}$. Then it follows from (\ref{fform2}) that 
\begin{equation}
\label{pp}
m_2(\alpha_2,\beta,T)=m_2(\alpha_2,\beta^{\ast},T).
\end{equation}
According to (\ref{mn01})---(\ref{mn03}), the sign of the function $m_2(\alpha_2,\beta,t)-\frac{t}{2}$, $t\in\mathbb{R}$, is the same as the sign of the number $1-\alpha_2^2-\beta^2$.
Then the numbers $1-\alpha_2^2-\beta^2$ and $1-\alpha_2^2-\beta^{\ast 2}$ have the same sign.
If these numbers are zero, then $|\beta|=|\beta^{\ast}|$.

If $1-\alpha_2^2-\beta^2>0$, then $|\beta|=|\beta^{\ast}|$ due to (\ref{pp}), (\ref{mn02}) and the increase of the function $\frac{\sinh x}{x}$ for $x>0$.

If $1-\alpha_2^2-\beta^2<0$, then $T<\frac{2\pi}{w_2}$ due to Proposition 9 and the condition $g\notin{\rm Cut}_2^{{\rm loc}}$. Then $|\beta|=|\beta^{\ast}|$ by virtue of (\ref{pp}), (\ref{mn03}) and the decrease of the function 
$\frac{\sin x}{x}$ on the interval $(0,\pi)$.

Let $\beta=\beta^{\ast}$. It follows from (\ref{fform2}) and (\ref{pp}) that either $\varphi_0^{\ast}=\varphi_0+2\pi n$, $n\in\mathbb{Z}$, or $B=0$. 
In the first case, the shortest arcs
$\tilde{\gamma}_2(\varphi_0,\alpha_2,\beta;t)$, $\tilde{\gamma}_2(\varphi_0^{\ast},\alpha_2^{\ast},\beta^{\ast};t)$,  $0\leq t\leq T$,
coincide, which is impossible. In the second case 
$\beta^2>1-\alpha_2^2$ and $T=\frac{2\pi}{w_2}$, i.e., 
$g\in{\rm Cut}_2^{{\rm loc}}$ by virtue of Corollary 1.

Let $\beta^{\ast}=-\beta$, $\beta>0$.  It follows from  (\ref{mn01})---(\ref{mn03}), and (\ref{form1})  that ${\rm Im}(A)$ is odd with respect to $\beta$. Therefore, ${\rm Im}(A)=0$.

If $\beta^2=1-\alpha_2^2$ then due to (\ref{mn01}),  the equality ${\rm Im}(A)=0$ 
can be rewritten in the form
$\tan\frac{\beta T}{2}=\frac{\beta T}{2}$. Then $\pi<\frac{\beta T}{2}<\frac{3\pi}{2}$ and, taking into account (\ref{form1}), we obtain
${\rm Re}(A)=\cos\frac{\beta T}{2}+\frac{\beta T}{2}\sin\frac{\beta T}{2}<0.$

If $\beta^2<1-\alpha_2^2$ then due to (\ref{mn02}), the equality ${\rm Im}(A)=0$ can be rewritten as
\begin{equation}
\label{pp2}
\tan kx=k{\rm tanh}x,\quad k=\frac{|\beta|}{\sqrt{1-\alpha_2^2-\beta^2}}>0,\,\, x=\frac{T\sqrt{1-\alpha_2^2-\beta^2}}{2}>0.
\end{equation}
Note that the function $f(x)=\tan kx-k{\rm th}x$ 
is increasing on
$(0,+\infty)$, $f(\pi/k)<0$,
$$\lim\limits_{x\rightarrow 0}f(x)=0,\,\, \lim\limits_{x\rightarrow \frac{\pi}{2k}-0}f(x)=+\infty,\,\,\lim\limits_{x\rightarrow \frac{\pi}{2k}+0}f(x)=-\infty,\,\,
\lim\limits_{x\rightarrow \frac{3\pi}{2k}-0}f(x)=+\infty.$$
Thus $f(x)>0$ for $x\in\left(0,\frac{\pi}{2k}\right)$, $f(x)<0$ for
$x\in\left(\frac{\pi}{2k},\frac{\pi}{k}\right]$ and $f(x)$ 
has the unique root in the interval
$\left(\frac{\pi}{k},\frac{3\pi}{2k}\right)$. Then $\pi<\frac{\beta T}{2}<\frac{3\pi}{2}$. 
Taking into account
(\ref{mn02}), (\ref{form1}) and the equality ${\rm Im}(A)=0$, 
\begin{equation}
\label{pp3}
\sin\frac{\beta T}{2}=\frac{-\beta m_2}{\sqrt{\beta^2m_2^2+n_2^2}},\,\,
\cos\frac{\beta T}{2}=\frac{-n_2}{\sqrt{\beta^2m_2^2+n_2^2}}\,\,\Rightarrow\,\, {\rm Re}(A)=-\sqrt{\beta^2m_2^2+n_2^2}<0.
\end{equation}

Let $\beta^2>1-\alpha_2^2$. If $n_2(T)=0$, then $T=\frac{\pi}{\sqrt{\beta^2+\alpha_2^2-1}}$, $m_2(T)=\frac{1}{\sqrt{\beta^2+\alpha_2^2-1}}$ by virtue of (\ref{mn03}). 
Then it follows from the equality ${\rm Im}(A)=0$ that 
$\cos\frac{\beta T}{2}=0$, i.e., $\frac{\beta}{\sqrt{\beta^2+\alpha_2^2-1}}=1+2l$, $l\in\mathbb{N}$. 

If $l=1$ then ${\rm Re}(A)=\beta m_2(T)\sin\frac{\beta T}{2}<0.$

If $l>1$ then $\frac{\beta}{\sqrt{1-\alpha_2^2}}<\frac{3}{2\sqrt{2}}$ and $T>\frac{3\pi}{\beta}$. 
From further reasoning it will follow that in the case the  geodesic segment $\tilde{\gamma}_2(\varphi_0,\alpha_2,\beta;t)$, $0\leq t\leq T$, is not the shortest arc.

If $n_2(T)\neq 0$ then by virtue of (\ref{mn03}) the equality ${\rm Im}(A)=0$ can be rewritten in the form
\begin{equation}
\label{pp0}
\tan kx=k\tan x,\quad k=\frac{|\beta|}{\sqrt{\beta^2+\alpha_2^2-1}}>1,\,\, x=\frac{T\sqrt{\beta^2+\alpha_2^2-1}}{2}\in \left(0,\pi\right).
\end{equation}

In the proof of Theorem~6 in [8] it is shown that the equation (\ref{pp0}) has no solutions for $k\in (1,2]$ and $k=3$. The smallest positive root of this equation belongs to the interval  $\left(\frac{3\pi}{2k},\frac{\pi}{k-1}\right)$ for $2<k<3$, $\left(\frac{\pi}{k-1},\frac{3\pi}{2k}\right)$ for $k>3$. Consequently, $\cos x<0$, $\cos kx>0$ for $2<k<3$ and $\cos x>0$, $\cos kx<0$ for $k>3$. 
Hence, taking into account (\ref{mn03}), (\ref{form1}) and the equality ${\rm Im}(A)=0$, we again obtain (\ref{pp3}).
Then $g\in{\rm Cut}_2^{{\rm glob}}$, which is impossible.
\noindent $\square$

\vspace{2mm}

{\bf Proposition 11.} {\sl 
	Let $\tilde{\gamma}_1(\varphi_0,\alpha_2,\beta;t)$, $0\leq t\leq T$, be a noncontinuable shortest arc
	in $(SU(1,1)\times\mathbb{R},d_1)$.
	Then $\tilde{\gamma}_2(\varphi_0,\alpha_2,\beta+\alpha_2;t)$, $0\leq t\leq T$, is a noncontinuable shortest arc in $(SU(1,1)\times\mathbb{R},d_2)$. The converse is also true.
}

\vspace{2mm}

{\sl Proof.} 
Denote $g_1=(A,B,v)=\tilde{\gamma}_1(\varphi_0,\alpha_2,\beta;T)$. By Theorem 4, $v=\alpha_2 T$. Then $g_2:=\tilde{\gamma}_2(\varphi_0,\alpha_2,\beta+\alpha_2;T)=\left(A e^{iv/2},Be^{-iv/2},v\right)$ by virtue of 
Proposition~2 and (\ref{Basis1}). Suppose that $\tilde{\gamma}_2(\varphi_0^{\ast},\alpha_2^{\ast},\beta^{\ast};t)$, $0\leq t\leq T_1<T$, is a shortest arc in $(SU(1,1)\times\mathbb{R},d_2)$ joining ${\rm Id}$ and $g_2$; $v=\alpha_2^{\ast}T_1$ by Theorem 4. By virtue of Proposition~2, the geodesic segment
$\tilde{\gamma}_1(\varphi_0^{\ast},\alpha_2^{\ast},\beta^{\ast}-\alpha_2^{\ast};t)$, $0\leq t\leq T_1$, in $(SU(1,1)\times\mathbb{R},d_1)$ joins ${\rm Id}$ and $g_1$, which is impossible. Hence $\tilde{\gamma}_2(\varphi_0,\alpha_2,\beta+\alpha_2;t)$, $0\leq t\leq T$, is a shortest arc in $(SU(1,1)\times\mathbb{R},d_2)$.  
It is easy to check that if the geodesic segment 
$\tilde{\gamma}_2(\varphi_0,\alpha_2,\beta+\alpha_2;t)$,  $0\leq t\leq T_2$, $T_2>T$, is a shortest arc in $(SU(1,1)\times\mathbb{R},d_2)$, then 
$\tilde{\gamma}_1(\varphi_0,\alpha_2,\beta;T)$, $0\leq t\leq T_2$, $T_2>T$, is a shortest arc in $(SU(1,1)\times\mathbb{R},d_1)$, which is impossible. Thus $\tilde{\gamma}_2(\varphi_0,\alpha_2,\beta+\alpha_2;t)$,  $0\leq t\leq T$, is a noncontinuable shortest arc in $(SU(1,1)\times\mathbb{R},d_2)$.

The converse statement is proved in a similar way.
\noindent $\square$

\vspace{2mm}

{\bf Proposition 12.} {\sl 
	The cut locus ${\rm Cut}_1$ of the sub-Riemannian space $(SU(1,1)\times\mathbb{R},d_1)$ corresponding to ${\rm Id}$ is the set
	$${\rm Cut}_1={\rm Cut}_1^{{\rm loc}}\cup{\rm Cut}_1^{{\rm glob}},$$
	where
	$${\rm Cut}_1^{{\rm loc}}=\{(A,B,v)\in SU(1,1)\times\mathbb{R}\mid A\neq e^{-iv/2},\,\,B=0,\,\,v\in\mathbb{R}\},$$
	$${\rm Cut}_1^{{\rm glob}}=\{(A,B,v)\in SU(1,1)\times\mathbb{R}\mid B\neq 0,\,\,\pi-v/2\in{\rm Arg}(A),\,\,v\in\mathbb{R}\},$$
	and ${\rm Arg}(A)$ is the set of all values of the argument of a complex number $A$.}

\vspace{2mm}

{\sl Proof.}
It follows from the definition of the cut set, Propositions~2 and 11, and Theorem~4 that
$$g_2=(A,B,v)\in{\rm Cut}_2\,\,\Leftrightarrow\,\, g_1=(A,B,v)\exp(-vE_3)=\left(Ae^{-iv/2},Be^{iv/2},v\right)\in {\rm Cut}_1.$$
From here and Corollary~1 it follows that $g_2\in{\rm Cut}_2^{{\rm loc}}$ if and only if 
$g_1\in {\rm Cut}_1^{{\rm loc}}$.

If $g_2\in{\rm Cut}_2^{{\rm glob}}$, i.e., ${\rm Im}(A)=0$ and ${\rm Re}(A)<-1$ by virtue of Proposition~10, then $Ae^{-iv/2}=|{\rm Re}(A)|e^{i(\pi-v/2)}$, $|B|^2={\rm Re}^2(A)-1>0$; therefore $g_1\in {\rm Cut}_1^{{\rm glob}}$.

Conversely, if $g_1:=(A_0,B_0,v_0)\in{\rm Cut}_1^{{\rm glob}}$, i.e., $B_0\neq 0$ and $\pi-\frac{v_0}{2}\in{\rm Arg}(A_0)$, then
$$B=B_0e^{-iv_0/2}\neq 0,\quad A=A_0e^{iv_0/2}=|A_0|e^{i(\pi-v_0/2)}e^{iv_0/2}=|A_0|e^{i\pi}=-|A_0|<-1,$$
Hence, $g_2=(A,B,v)\in {\rm Cut}_1^{{\rm glob}}$.
\noindent $\square$

\section{Noncontinuable shortest arcs on the Lie group $SU(1,1)\times\mathbb{R}$}

By Proposition 4, to find shortest arcs in $(SU(1,1)\times\mathbb{R},\rho_i)$, $i=1,2$, it suffices study the geodesic segments $\tilde{\gamma}_i(t)$, $0\leq t\leq T$, such that $\tilde{\gamma}_i(0)={\rm Id}$.

\vspace{2mm}

{\bf Theorem 5.} {\sl
	Let $\alpha_2\neq\pm 1$, $\beta\neq 0$ and $\tilde{\gamma}_2(t)=\tilde{\gamma}_2(\varphi_0,\alpha_2,\beta;t)$, $0\leq t\leq T$, is a noncontinuable shortest arc of the sub--Riemannian space $(SU(1,1)\times\mathbb{R},d_2)$. Then
	
	1. If $\frac{|\beta|}{\sqrt{1-\alpha_2^2}}\geq\frac{2}{\sqrt{3}}$ then $T=\frac{2\pi}{\sqrt{\beta^2+\alpha_2^2-1}}$.
	
	2. If $\frac{|\beta|}{\sqrt{1-\alpha_2^2}}=1$ then $T\in \left(\frac{2\pi}{|\beta|},\frac{3\pi}{|\beta|}\right)$ and $T$ 
	satisfies the system of equations
	$$\cos\frac{\beta T}{2}=\frac{-2}{\sqrt{4+\beta^2T^2}},\quad \sin\frac{\beta T}{2}=\frac{-\beta T}{\sqrt{4+\beta^2T^2}}.$$
	
	3. If $0<\frac{|\beta|}{\sqrt{1-\alpha_2^2}}<1$ then $T\in\left(\frac{2\pi}{|\beta|},\frac{3\pi}{|\beta|}\right)$ and  
	$T$ satisfies the system of equations
	$$\cos kx=\frac{-1}{\sqrt{1+k^2\tanh^2x}},\quad \sin kx=\frac{-k \tanh x}{\sqrt{1+k^2\tanh^2x}},$$
	where $k$ and $x$  are defined by formulas (\ref{pp2}).
	
	4. If $\frac{|\beta|}{\sqrt{1-\alpha_2^2}}=\frac{3}{2\sqrt{2}}$ then $T=\frac{3\pi}{|\beta|}$.
	
	5. If $\frac{3}{2\sqrt{2}}<\frac{|\beta|}{\sqrt{1-\alpha_2^2}}<\frac{2}{\sqrt{3}}$ then
	$\frac{3\pi}{|\beta|}< T<\frac{2\pi\left(|\beta|+\sqrt{\beta^2+\alpha_2^2-1}\right)}{\sqrt{1-\alpha_2^2}}<\frac{4\pi}{|\beta|}$
	and  $T$ satisfies the system of equations
	$$\cos kx=\frac{1}{\sqrt{1+k^2\tan^2x}},\quad
	\sin kx=\frac{k\tan x}{\sqrt{1+k^2\tan^2 x}}<0,$$
	where $k$ and $x$  are defined by formulas (\ref{pp0}).
	
	6. If $1<\frac{|\beta|}{\sqrt{1-\alpha_2^2}}<\frac{3}{2\sqrt{2}}$ then
	$\frac{2\pi}{\mid\beta\mid}<\frac{2\pi\left(|\beta|+\sqrt{\beta^2+\alpha_2^2-1}\right)}{\sqrt{1-\alpha_2^2}}<T<\frac{3\pi}{|\beta|}$
	and  $T$ satisfies the system of equations
	$$\cos kx=\frac{-1}{\sqrt{1+k^2\tan^2x}},\quad
	\sin kx=\frac{-k\tan x}{\sqrt{1+k^2\tan^2 x}}<0,$$
	where $k$ and $x$  are defined by formulas (\ref{pp0}).
}

\vspace{2mm}

{\sl Proof.}
Let $\alpha_2\neq\pm 1$, $\beta\neq 0$ and $\tilde{\gamma}_2(t)=\gamma_2(\varphi_0,\alpha_2,\beta;t)$,  $0\leq t\leq T$,  be a noncontinuable shortest arc. By virtue of Proposition~10, $g=(A,B,v):=\tilde{\gamma}_2(T)$
belongs to the union of ${\rm Cut}_2^{{\rm loc}}$ and ${\rm Cut}_2^{{\rm glob}}$.

If $g\in {\rm Cut}_2^{{\rm loc}}$ then $\beta^2>1-\alpha_2^2$, $m_2(T)=0$ and $T=\frac{2\pi}{\sqrt{\beta^2+\alpha_2^2-1}}$ by Theorem~3.

Let now $g\in {\rm Cut}_2^{{\rm glob}}$. Then by virtue of (\ref{form1}),
\begin{equation}
\label{qq}
n_2\sin\frac{\beta T}{2}-\beta m_2\cos\frac{\beta T}{2}=0,\quad
n_2\cos\frac{\beta T}{2}+\beta m_2\sin\frac{\beta T}{2}<0.
\end{equation}

If $\beta^2=1-\alpha_2^2$ then, taking into account (\ref{mn01}), the conditions (\ref{qq})  can be rewritten as $\tan\frac{\beta T}{2}=\frac{\beta T}{2}$, $\cos\frac{\beta T}{2}<0$,  whence follows $\pi<\frac{\beta T}{2}<\frac{3\pi}{2}$. Then
$$\cos\frac{\beta T}{2}=\frac{-1}{\sqrt{1+\tan^2(\beta T/2)}}=\frac{-2}{\sqrt{4+\beta^2T^2}},\quad \sin\frac{\beta T}{2}=\frac{-\beta T}{\sqrt{4+\beta^2T^2}},$$
and item~2 of Theorem~5 is proved.

If $\beta^2<1-\alpha_2^2$ then, taking into account  (\ref{mn02}), the conditions (\ref{qq}) can be rewritten as (\ref{pp2}), and besides $\cos kx<0$. 
In the proof of Proposition~10 it is shown that then $\pi<kx=\frac{|\beta| T}{2} <\frac{3\pi}{2}$, therefore
$$\cos kx=\frac{-1}{\sqrt{1+\tan^2kx}}=\frac{-1}{\sqrt{1+k^2\tanh^2x}},\quad\sin kx=\frac{-k{\rm tanh}x}{\sqrt{1+k^2\tanh^2x}},$$
and item~3 of Theorem~5 is proved.

Let $\beta^2>1-\alpha_2^2$. If $n_2(T)=0$ then, as shown in the proof of Proposition~10,
$T=\frac{\pi}{\sqrt{\beta^2+\alpha_2^2-1}}$ and $\frac{|\beta|}{\sqrt{\beta^2+\alpha_2^2-1}}=3$.
This implies item~4 of Theorem~5.

Let $n_2(T)\neq 0$. Taking into account (\ref{mn03}), the conditions (\ref{qq}) can be rewritten as (\ref{pp0}), and besides $\cos kx<0$. In the proof of Theorem~6 in [8] it is shown that
the equation (\ref{pp0}) has no solutions for $1<k\leq 2$ and $k=3$. 

If  $2<k<3$ then the smallest positive root of the equation (\ref{pp0}) then the smallest positive root of the equation belongs to the interval
$\left(\frac{3\pi}{2k},\frac{\pi}{k-1}\right)$. Taking into account (\ref{pp0}), it means that
$$\frac{3}{2\sqrt{2}}<\frac{|\beta|}{\sqrt{1-\alpha_2^2}}<\frac{2}{\sqrt{3}},\quad \frac{3\pi}{|\beta|}<T<\frac{2\pi\left(|\beta|+\sqrt{\beta^2+\alpha_2^2-1}\right)}{\sqrt{1-\alpha_2^2}}<\frac{4\pi}{|\beta|}.$$
Then $\cos x<0$, $\tan x<0$, $\cos kx>0$ and by the virtue of (\ref{qq}),
$$\cos kx=\frac{1}{\sqrt{1+\tan^2 kx}}=\frac{1}{\sqrt{1+k^2\tan^2 x}},\quad
\sin kx=\frac{k\tan x}{\sqrt{1+k^2\tan^2 x}}.$$
Item~5 of Theorem~5 is proved.

If $k>3$, then the smallest positive root of the equation (\ref{pp0}) belongs to the interval $\left(\frac{\pi}{k-1},\frac{3\pi}{2k}\right)$. Taking into account (\ref{pp0}), it means that
$$1<\frac{|\beta|}{\sqrt{1-\alpha_2^2}}<\frac{3}{2\sqrt{2}},\quad
\frac{2\pi}{\mid\beta\mid}<\frac{2\pi\left(|\beta|+\sqrt{\beta^2+\alpha_2^2-1}\right)}{\sqrt{1-\alpha_2^2}}<T<\frac{3\pi}{|\beta|}.$$
Then $\cos x>0$, $\cos kx<0$ and by the virtue of (\ref{qq})
$$\cos kx=\frac{-1}{\sqrt{1+\tan^2 kx}}=\frac{-1}{\sqrt{1+k^2\tan^2 x}},\quad
\sin kx=\frac{-k\tan x}{\sqrt{1+k^2\tan^2 x}}.$$
Item~6 of Theorem~5 is proved.
\noindent $\square$

\vspace{2mm}

It follows from Proposition~11  that if $\tilde{\gamma}_1(\varphi_0,\alpha,\beta,t)$, with $0\leq t\leq T$,
is a noncontinuable shortest arc in $(SU(1,1)\times\mathbb{R},d_1)$, then $T$ satisfies items 1--6 of Theorem~5
after replacing $\beta$ by $\beta+\alpha_2$.

\vspace{1mm}

The following  is immediate from Propositions 2, 11.

\vspace{2mm}

{\bf Proposition 13.} {\sl 
For every $(A,B,v)\in SU(1,1)\times\mathbb{R}$,	
$$d_1({\rm Id},(A,B,v))=d_2({\rm Id},(Ae^{iv/2},Be^{-iv/2},v)).$$ 
}

The proof of the following theorem repeats the proof of Theorem 7 in [8].
We omit it.

\vspace{1mm}

{\bf Theorem 6.} {\sl 
	Let $\beta\neq 0$ and $\tilde{\gamma}_2(t)=\tilde{\gamma}_2(\varphi_0,\alpha_2,\beta;t)$,  $0\leq t\leq T$, be a noncontinuable shortest arc of the sub-Riemannian space $(SU(1,1)\times\mathbb{R},d_2)$. Then
	
	1) The function $T=T(|\beta|)$ decreases strictly on intervals
	$\left(0,\frac{3}{2\sqrt{2}}\sqrt{1-\alpha_2^2}\right]$, 
		\linebreak $\left[\frac{2}{\sqrt{3}}\sqrt{1-\alpha_2^2},+\infty\right)$ and increases strictly on
	$\left[\frac{3}{2\sqrt{2}}\sqrt{1-\alpha_2^2},\frac{2}{\sqrt{3}}\sqrt{1-\alpha_2^2}\right]$.
	
	2) The function $T=T(|\beta|)$ is continuous, piecewise real-analytic and
	$T(0,+\infty)=(0,+\infty)$.
	
	3) The function $T=T(|\beta|)$ has the local minimum $2\sqrt{2}\pi/\sqrt{1-\alpha_2^2}$ at
	$\frac{3}{2\sqrt{2}}\sqrt{1-\alpha_2^2}$  and the local maximum $2\sqrt{3}\pi/\sqrt{1-\alpha_2^2}$ at $\frac{2}{\sqrt{3}}\sqrt{1-\alpha_2^2}$.
}

\vspace{2mm}

The following theorem is proved in paper [12].

\vspace{2mm}

{\bf Theorem 7.} {\sl
	Let $\mathbb{H}\times\mathbb{R}$ be a connected $(n+1)$--dimensional Lie group p (with the unit ${\rm Id}=(e,1)$ and the Lie algebra $\mathfrak{g}\oplus\mathfrak{g}_1$) with the left-invariant sub--Riemannian metric
	$d$ generated by a completely nonholonomic distribution $D$ with
	$$D({\rm Id})={\rm span}(e_1,\dots,e_{n-1},e_{n+1})\subset\mathfrak{g}\oplus\mathfrak{g}_1,\quad e_1,\dots,e_{n-1}\in\mathfrak{g},\,\,e_{n+1}=1\in\mathfrak{g}_1,$$
	and the inner product $\langle\cdot,\cdot\rangle$ with the orthonormal basis $e_1,\dots,e_{n-1}$, $e_{n+1}$ 
	on $D({\rm Id})$. Then
	$d^2({\rm Id},(g,e^v))=v^2+d^2({\rm Id},(g,1))$ for all
	$(g,e^v)\in\mathbb{H}\times\mathbb{R}$. 
}

\vspace{2mm}

The following is immediate from Theorem 7 and (\ref{abc2}).

\vspace{2mm}

{\bf Proposition 14.} {\sl
	For every $(A,B,v)\in SU(1,1)\times\mathbb{R}$,
	$$d_2^2({\rm Id},(A,B,v))=v^2+d_2^2({\rm Id},(A,B,0)).$$ 
}

It follows from Remark 1 that 
$d_2({\rm Id},(A,B,v))=d_2({\rm Id},\tilde{\Omega}(A,B,v))$, where the left-invariant inner metric $d_2$ 
on the right side of the equality is defined on $SL(2,\mathbb{R})\times\mathbb{R}$ by Theorem~3.
In [9] the exact formulas for $d_2({\rm Id},\tilde{\Omega}(A,B,0))$ are found. Using them and Propositions 13, 14, we can obtain exact formulas for $d_i({\rm Id},(A,B,v))$, $i=1,2$. We will not
present them because of their bulkiness.

\section{Sub-Riemannian geodesics on $SO_0(2,1)\times\mathbb{R}$}

The pseudoeuclidean space $\mathbb{E}^{2,1}$ is the vector space $\mathbb{R}^{3}$
with pseudoscalar product $\{x,y\}:=-x_1y_1+x_2y_2+x_3y_3$, where $x=(x_1,x_2,x_3)$, $y=(y_1,y_2,y_3)$. 
The Lorentz group $SO_0(2,1)$ is the connected component
of the unit in the group  $P(2,1)$ 
of all linear preserving the pseudoscalar product $\{\cdot,\cdot\}$ 
transformations of the space $\mathbb{E}^{2,1}$.
It consists of those elements in $P(2,1)$  which simultaneously preserve the time direction (the first coordinate) and orientation of the space $\mathbb{E}^{2,1}$. Moreover,
$$SO_0(2,1)=\left\{g\in GL(3,\mathbb{R})\left|\,g^{T}Jg=J,\,\,J=\rm{diag}(-1,1,1),\,\,\rm{det}(g)=1,\,\,g_{11}>0\right.\right\}.$$
Its Lie algebra is
$\quad\mathfrak{so}(2,1)=\left\{A\in\mathfrak{gl}(3,\mathbb{R})\left|\,A^{T}J+JA=0,\,\,J={\rm diag}(-1,1,1)\right.\right\}.$

We will deal with the trivial abelian extension $SO_0(2,1)$:
$$SO_0(2,1)\times\mathbb{R}=\scriptsize{\left\{(C,v):=\left(\begin{array}{ccc}
	C &  & 0 \\
	&  & 0 \\
	0 & 0 & e^v 
	\end{array}\right)\,\,\left|\,\,C\in SO_0(2,1),\,\,v\in\mathbb{R}\right.\right\}.}$$
The Lie algebra  $\mathfrak{so}(2,1)\oplus\mathbb{R}$ of the Lie group
$SO_0(2,1)\times\mathbb{R}$ has the basis $E_1,E_2,E_3,E_4$ satisfying (\ref{abc}):
\begin{equation}
\label{basis2}
E_1=e_{13}+e_{31},\quad E_2=e_{12}+e_{21},\quad E_3=e_{23}-e_{32},\quad E_4=e_{44},
\end{equation}
where $e_{ij}$, with $i,j=1,\dots,4$, denotes the $4\times$--matrix,
having $1$ in the $i$th row and the $j$th column, while
all other entries are $0$.

For convenience we denote by $\rho_i$, $i=1,2$, the left-invariant sub-Riemannian metric on $SO_0(2,1)\times\mathbb{R}$ ((with the unit ${\rm Id}$) defined by the completely nonholonomic left-invariant distribution $\Delta_i$ with $\Delta_i({\rm Id})={\rm span}(e_1,e_2,e_3)$, where the vectors $e_1,e_2,e_3$ are given by formulas (\ref{abc1}) for $i=1$ and (\ref{abc2}) for $i=2$; 
we define the inner product $\langle\cdot,\cdot\rangle_i$ on $\Delta_i({\rm Id})$ with the orthonormal basis
$e_1,e_2,e_3$.

\vspace{2mm}

{\bf Theorem 8.} {\sl
	Let
	$$\mu_2=t,\quad \nu_2=t^2/2,\quad\mbox{if }\,\,\beta^2=1-\alpha_2^2;$$
	$$\mu_2=\frac{\sinh\left(t\sqrt{1-\alpha_2^2-\beta^2}\right)}{\sqrt{1-\alpha_2^2-\beta^2}},\quad \nu_2=\frac{\cosh\left(t\sqrt{1-\alpha_2^2-\beta^2}\right)-1}{1-\alpha_2^2-\beta^2},\quad\mbox{if }\,\,\beta^2<1-\alpha_2^2;$$
	$$\mu_2=\frac{\sin\left(t\sqrt{\beta^2+\alpha_2^2-1}\right)}{\sqrt{\beta^2+\alpha_2^2-1}},\quad \nu_2=\frac{1-\cos\left(t\sqrt{\beta^2+\alpha_2^2-1}\right)}{\beta^2+\alpha_2^2-1},\quad\mbox{if }\,\,\beta^2>1-\alpha_2^2.$$
	If $\alpha_2\neq\pm 1$ then the geodesic $\gamma_2(\alpha_1,\alpha_2,\alpha_3,\beta;t)$, 
	$t\in\mathbb{R}$, of the sub-Riemannian space $(SO_0(2,1)\times\mathbb{R},\rho_2)$ defined by the formula (\ref{geod2}), is equal to $(C,v)(t)$, where $v(t)=\alpha_2t$ and the columns $C_j(t)$, $j=1,2,3$, of the
	matrix $C(t)\in SO_0(2,1)$ have the form
	$$C_1(t)=\left(
	\begin{array}{c}
	1+\nu_2(1-\alpha_2^2)  \\
	\sqrt{1-\alpha_2^2}\left(\mu_2\cos{\varphi_0}-\beta\nu_2\sin{\varphi_0}\right)  \\
	\sqrt{1-\alpha_2^2}\left(\mu_2\sin{\varphi_0}+\beta\nu_2\cos{\varphi_0}\right)
	\end{array}
	\right),$$
	$$C_2(t)=\left(
	\begin{array}{c}
	\sqrt{1-\alpha_2^2}\left(\mu_2\cos(\beta t+\varphi_0)+\beta \nu_2\sin(\beta t+\varphi_0)\right)  \\
	(1-\beta^2\nu_2)\cos\beta t+\beta \mu_2\sin\beta t+(1-\alpha_2^2)\nu_2\cos(\beta t+\varphi_0)\cos{\varphi_0}  \\
	-(1-\beta^2\nu_2)\sin\beta t+\beta \mu_2\cos\beta t+(1-\alpha_2^2)\nu_2\cos(\beta t+\varphi_0)\sin{\varphi_0}
	\end{array}
	\right),$$
	$$C_3(t)=\left(
	\begin{array}{c}
	\sqrt{1-\alpha_2^2}\left(\mu_2\sin(\beta t+\varphi_0)-\beta \nu_2\cos(\beta t+\varphi_0)\right)  \\
	(1-\beta^2 \nu_2)\sin\beta t-\beta \mu_2\cos\beta t+(1-\alpha_2^2)\nu_2\sin(\beta t+\varphi_0)\cos{\varphi_0}  \\
	(1-\beta^2 \nu_2)\cos\beta t+\beta \mu_2\sin\beta t+(1-\alpha_2^2)\nu_2\sin(\beta t+\varphi_0)\sin{\varphi_0}
	\end{array}
	\right),$$
	\begin{equation}
	\label{ss}
	\alpha_1=\sqrt{1-\alpha_2^2}\sin\varphi_0,\quad \alpha_3=\sqrt{1-\alpha_2^2}\cos\varphi_0.
	\end{equation}
	
	If $\alpha_2=\pm 1$ the the geodesic of the sub-Riemannian space $(SO_0(2,1)\times\mathbb{R},\rho_2)$	
	defined by the formula (\ref{geod2}),  is equal to $\gamma_2(t)=(E,v)(t)$, $t\in\mathbb{R}$, where $E$ is the identity $3\times 3$--matrix, $v(t)=\alpha_2t$.	
}

\vspace{1mm}

{\sl Proof.}
Let $\alpha_2\neq\pm 1$. Using (\ref{norm}), (\ref{abc2}), (\ref{basis2}), it is easy to get that
$$\exp(t(\alpha_1e_1+\alpha_2e_2+\alpha_3e_3+\beta e_4))
=\left(\begin{array}{cccc}
1+(1-\alpha_2^2)\nu_2 & \alpha_3\mu_2+\alpha_1\beta\nu_2 & \alpha_1\mu_2-\alpha_3\beta\nu_2 & 0  \\
\alpha_3\mu_2-\alpha_1\beta\nu_2 & 1+(\alpha_3^2-\beta^2)\nu_2 & \alpha_1\alpha_3\nu_2-\beta\mu_2 & 0  \\
\alpha_1\mu_2+\alpha_3\beta\nu_2 & \alpha_1\alpha_3\nu_2+\beta\mu_2 & 1+(\alpha_1^2-\beta^2)\nu_2 & 0 \\
0 & 0 & 0 & \alpha_2t
\end{array}
\right),$$
$$\exp(-t\beta e_4)=
\left(\begin{array}{cccc}
1 & 0 & 0 & 0  \\
0 & \cos\beta t & \sin\beta t & 0  \\
0 & -\sin\beta t & \cos\beta t & 0 \\
0 & 0 & 0 & 1
\end{array}
\right).$$
By (\ref{geod2}),  it remains to multiply these matrices.

In the case $\alpha_2=\pm 1$ the statement of Theorem 8  follows from Proposition~3, (\ref{abc2}), (\ref{basis2}).
\noindent $\square$

\vspace{2mm}

The following is immediate from Theorem 8 and Propositions 2, 3.

\vspace{2mm}

{\bf Theorem 9.} {\sl
	Let
	$$\mu_1=t,\quad \nu_1=t^2/2,\quad\mbox{if }(\beta+\alpha_2)^2=1-\alpha_2^2;$$
	$$\mu_1=\frac{\sinh\left(t\sqrt{1-\alpha_2^2-(\beta+\alpha_2)^2}\right)}{\sqrt{1-\alpha_2^2-(\beta+\alpha_2)^2}},\quad \nu_1=\frac{\cosh\left(t\sqrt{1-\alpha_2^2-(\beta+\alpha_2)^2}\right)-1}{1-\alpha_2^2-(\beta+\alpha_2)^2},$$
	if $(\beta+\alpha_2)^2<1-\alpha_2^2;$
	$$\mu_1=\frac{\sin\left(t\sqrt{(\beta+\alpha_2)^2+\alpha_2^2-1}\right)}{\sqrt{(\beta+\alpha_2)^2+\alpha_2^2-1}},\quad\nu_1=\frac{1-\cos\left(t\sqrt{(\beta+\alpha_2)^2+\alpha_2^2-1}\right)}{(\beta+\alpha_2)^2+\alpha_2^2-1},$$
	if $(\beta+\alpha_2)^2>1-\alpha_2^2.$
	
	If $\alpha_2\neq\pm 1$ then  the geodesic $\gamma_1(\alpha_1,\alpha_2,\alpha_3,\beta;t)$,
	$t\in\mathbb{R}$, of the sub-Riemannian space
	$(SO_0(2,1)\times\mathbb{R},\rho_1)$ defined by the formula (\ref{geod}) is equal to
	$(C,v)(t)$, $t\in\mathbb{R}$, where $v(t)=\alpha_2t$ and   the columns $C_j(t)$, with $j=1,2,3$, of the
	matrix $C(t)\in SO_0(2,1)$ have the form	
	$$C_1(t)=\left(
	\begin{array}{c}
	1+\nu_2(1-\alpha_2^2)  \\
	\sqrt{1-\alpha_2^2}\left(\mu_2\cos{\varphi_0}-\beta\nu_2\sin{\varphi_0}\right)  \\
	\sqrt{1-\alpha_2^2}\left(\mu_2\sin{\varphi_0}+\beta\nu_2\cos{\varphi_0}\right)
	\end{array}
	\right),$$
	$$C_2(t)=\left(
	\begin{array}{c}
	\sqrt{1-\alpha_2^2}\left(\mu_2\cos(\beta t+\varphi_0)+(\beta+\alpha_2) \nu_2\sin(\beta t+\varphi_0)\right)  \\
	(1-(\beta+\alpha_2)^2\nu_2)\cos\beta t+(\beta+\alpha_2) \mu_2\sin\beta t+(1-\alpha_2^2)\nu_2\cos(\beta t+\varphi_0)\cos{\varphi_0}  \\
	-(1-(\beta+\alpha_2)^2\nu_2)\sin\beta t+(\beta+\alpha_2) \mu_2\cos\beta t+(1-\alpha_2^2)\nu_2\cos(\beta t+\varphi_0)\sin{\varphi_0}
	\end{array}
	\right),$$
	$$C_3(t)=\left(
	\begin{array}{c}
	\sqrt{1-\alpha_2^2}\left(\mu_2\sin(\beta t+\varphi_0)-(\beta+\alpha_2) \nu_2\cos(\beta t+\varphi_0)\right)  \\
	(1-(\beta+\alpha_2)^2 \nu_2)\sin\beta t-(\beta+\alpha_2) \mu_2\cos\beta t+(1-\alpha_2^2)\nu_2\sin(\beta t+\varphi_0)\cos{\varphi_0}  \\
	(1-(\beta+\alpha_2)^2 \nu_2)\cos\beta t+(\beta+\alpha_2) \mu_2\sin\beta t+(1-\alpha_2^2)\nu_2\sin(\beta t+\varphi_0)\sin{\varphi_0}
	\end{array}
	\right),$$
	and (\ref{ss}) holds.
	
	If $\alpha_2=\pm 1$ then the geodesic of the sub-Riemannian space $(SO_0(2,1)\times\mathbb{R},\rho_1)$
	defined by the formula  (\ref{geod}) is equal to $\gamma_1(t)=(C,v)(t)$,  $t\in\mathbb{R}$, where $v(t)=\alpha_2t$ and
	$$C(t)=\left(
	\begin{array}{ccc}
	1 & 0 & 0  \\
	0 & \cos{\alpha_2 t} & -\sin{\alpha_2 t}  \\
	0 & \sin{\alpha_2 t} & \cos{\alpha_2 t} 
	\end{array}
	\right).$$
}

\vspace{2mm}

{\bf Remark 2.} {\sl
	In the case $\alpha_2\neq\pm 1$ we will also denote the geodesic $\gamma_i(\alpha_1,\alpha_2,\alpha_3,\beta;t)$ by $\gamma_i(\varphi_0,\alpha_2,\beta;t)$,
	if (\ref{ss}) holds.
}

\section{The Cut Locus and First Conjugate Locus in $SO_0(2,1)\times\mathbb{R}$}

The Lie group $SU(1,1)$ is a simply connected covering of the group $SO_0(2,1)$; the double covering $\Pi:SU(1,1)\rightarrow SO_0(2,1)$  can be given as follows:
\begin{equation}
\label{formula}
\Pi(A,B)=\left(
\begin{array}{ccc}
A_1^2+A_2^2+B_1^2+B_2^2 & 2(A_1B_2-A_2B_1) & 2(A_1B_1+A_2B_2)  \\
2(A_1B_2+A_2B_1) & A_1^2-A_2^2-B_1^2+B_2^2 & 2(A_1A_2+B_1B_2) \\
2(A_1B_1-A_2B_2) & 2(B_1B_2-A_1A_2) & A_1^2-A_2^2+B_1^2-B_2^2
\end{array}
\right),
\end{equation}
where $A_1,A_2,B_1$, and $B_2$ are defined by (\ref{AB}). Then the mapping
$$\tilde{\Pi}:SU(1,1)\times\mathbb{R}\rightarrow SO_0(2,1)\times\mathbb{R},\quad \tilde{\Pi}(A,B,v)=(\Pi(A,B),v),$$
is a double covering of the Lie group $SO_0(2,1)\times\mathbb{R}$ 
by $SU(1,1)\times\mathbb{R}$.

It is easy to check that the differential $d\tilde{\Pi}({\rm Id})$
is the isomorphism of the Lie algebras $\mathfrak{su}(1,1)\oplus\mathbb{R}$
and $\mathfrak{so}(2,1)\oplus\mathbb{R}$ translating the basis $E_1,E_2,E_3,E_4$ of the Lie algebra $\mathfrak{su}(1,1)\oplus\mathbb{R}$ given by (\ref{Basis1}) to the same name basis of $\mathfrak{so}(2,1)\oplus\mathbb{R}$ given by (\ref{basis2}). Therefore, according to Theorems  1,  2, the mapping $\tilde{\Pi}:(SU(1,1)\times\mathbb{R},d_i)\rightarrow (SO_0(2,1)\times\mathbb{R},\rho_i)$, with $i=1,2$,  is a local isometry, and for all $(\alpha_1,\alpha_2,\alpha_3)\in\mathbb{S}^2$, $\beta,t\in\mathbb{R}$ we have
\begin{equation}
\label{soot}
\tilde{\Pi}(\tilde{\gamma}_i(\alpha_1,\alpha_2,\alpha_3,\beta;t))
=\gamma_i(\alpha_1,\alpha_2,\alpha_3,\beta;t),\quad i=1,2.
\end{equation}

Equality (\ref{formula}) and Corollary 1 imply

\vspace{2mm}

{\bf Corollary 2.} {\sl
	1.  The $n$-th conjugate time $t^{n}_{{\rm conj}}$ along the geodesic $\gamma_i(\alpha_1,\alpha_2,\alpha_3,\beta;t)$, where $\alpha_2\neq \pm 1$ and $w_i>0$ ($w_i$ are given by formula (\ref{w12})) of the sub-Riemannian space
	$(SO_0(2,1)\times\mathbb{R},\rho_i)$, $i=1,2$, has the form
	$$t^{2m-1}_{{\rm conj}}=\frac{2\pi m}{w_i},\quad t^{2m}_{{\rm conj}}=\frac{2x_m}{w_i},\quad m\in\mathbb{N},$$
	where $\{x_1,x_2,\dots\}$ are the ascending positive roots of the equation $\tan{x}=x$.
	
	2.  The first conjugate locus ${\rm Conj}_i^{1}$ of the sub-Riemannian space  $(SO_0(2,1)\times\mathbb{R},\rho_i)$, $i=1,2$, has the form
	$${\rm Conj}_1^{1}={\rm Conj}_2^{1}=\left\{(C,v)\in SO_0(2,1)\times\mathbb{R}\left|\,\,
	C=\left(
	\begin{array}{ccc}
	1 & 0 & 0 \\
	0 & \cos\psi & -\sin\psi \\
	0 & \sin\psi & \cos\psi \\
	\end{array}
	\right),\,\,\psi,\,v\in\mathbb{R}\right.\right\}.$$
}

\vspace{1mm}

The proof of the following is similar to the proof of Proposition 9 in [12]; we will omit it.

\vspace{2mm}

{\bf Proposition 15.} {\sl
	Let $\gamma_i(t)$, $0\leq t\leq T$, be a noncontinuable shortest arc 
	in  $(SO_0(2,1)\times\mathbb{R},\rho_i)$, $i=1,2$, such that $\gamma_i(0)={\rm Id}$. There exists only one shortest arc
	 $\tilde{\gamma}_i(t)$, $0\leq t\leq T$, in $(SU(1,1),d_i)$ such that  $\tilde{\gamma}_i(0)={\rm Id}$ and 
	 $\tilde{\Pi}(\tilde{\gamma}_i(t))=\gamma_i(t)$ for every $t\in [0,T]$.
}

\vspace{2mm}

The following is immediate from (\ref{soot}) and Propositions 6, 7.

\vspace{2mm}

{\bf Proposition 16.} {\sl
	1. If $\alpha_2=\pm 1$ then the geodesic $\gamma_i(t)=\exp(t\alpha_2 e_2)$,  $t\in\mathbb{R}$, in  $(SO_0(2,1)\times\mathbb{R},\rho_i)$, $i=1,2$, is a metric line.
	
	2. The geodesic $\gamma_1(t)=\gamma_1(\varphi_0,\alpha_2,-\alpha_2;t)$,  $t\in\mathbb{R}$,  in $(SO_0(2,1)\times\mathbb{R},\rho_1)$ is a metric line.
	
	3. The geodesic $\gamma_2(t)=\gamma_2(\varphi_0,\alpha_2,0;t)$,  $t\in\mathbb{R}$, 
	in $(SO_0(2,1)\times\mathbb{R},\rho_2)$  is a metric line.
}

\vspace{2mm}

{\bf Proposition 17.} {\sl
	The cut locus ${\rm Cut}_2$ of the sub-Riemannian space $(SO_0(2,1)\times\mathbb{R},\rho_2)$, corresponding to ${\rm Id}$, is the set
	$${\rm Cut}_2={\rm Cut}_2^{{\rm loc}}\cup {\rm Cut}_2^{{\rm glob}},$$
	where
	$${\rm Cut}_2^{{\rm loc}}=\{\tilde{\Pi}(A,0,v)\mid (A,0,v)\in SU(1,1)\times\mathbb{R},\,\,A\neq \pm 1,\,\,v\in\mathbb{R}\}$$
	$$=\left\{(C,v)\in SO_0(2,1)\times\mathbb{R}\left|\,\,
	C=\left(
	\begin{array}{ccc}
	1 & 0 & 0 \\
	0 & \cos\psi & -\sin\psi \\
	0 & \sin\psi & \cos\psi \\
	\end{array}
	\right),\,\,\psi\in (0,2\pi),\,v\in\mathbb{R}\right.\right\},$$
	$${\rm Cut}_2^{{\rm glob}}=\{\tilde{\Pi}(A,B,v)\mid (A,0,v)\in SU(1,1)\times\mathbb{R},\,\,{\rm Re}(A)=0,\,\,B\neq 0,\,\,v\in\mathbb{R}\}$$
	$$=\left\{(C,v)\in SO_0(2,1)\times\mathbb{R}\left|\,\,c_{21}=-c_{12}, \,\, c_{31}=-c_{13},\,\,c_{23}=c_{32},\,\,c_{22}+c_{33}<-2,\,\,v\in\mathbb{R}\right.\right\}.$$	}

{\sl Proof.}
Let us prove that ${\rm Cut}_2^{{\rm loc}}\subset{\rm Cut}_2$.  Let $g\in {\rm Cut}_2^{{\rm loc}}$ and $\gamma_2(t)$,  $0\leq t\leq T$, be a shortest arc
in  $(SO_0(2,1)\times\mathbb{R},\rho_2)$ joining ${\rm Id}$ and $g$.
By virtue of Proposition 15, there exists only one shortest arc  $\tilde{\gamma}_2(t)$,  $0\leq t\leq T$, of the sub-Riemannian space $(SU(1,1)\times\mathbb{R},d_2)$ such that $\tilde{\gamma}_2(0)={\rm Id}$ and $\tilde{\Pi}(\tilde{\gamma}_2(t))=\gamma_2(t)$ for every $t\in [0,T]$. 
It follows from Corollary 2 and the proof of Proposition 10 that $\tilde{g}:=\tilde{\gamma}_2(T)$  belongs to
the cut locus (for ${\rm Id}$) in $(SU(1,1)\times\mathbb{R},d_2)$, 
and there is another shortest arc 
$\tilde{\gamma}^{\ast}_2(t)$, $0\leq t\leq T$, joining ${\rm Id}$ and $\tilde{g}$ such that $\tilde{\Pi}(\tilde{g})=g$. Thus the geodesic segment $\tilde{\Pi}(\tilde{\gamma}^{\ast}_2(t))\neq\gamma_2(t)$, $0\leq t\leq T$, in $(SO_0(2,1)\times\mathbb{R},\rho_2)$ is
a shortest arc joining ${\rm Id}$ and $g$. Then $g\in {\rm Cut}_2$ by Proposition 8.

Let us prove that ${\rm Cut}_2^{{\rm glob}}\subset {\rm Cut}_2$.  Let $g\in {\rm Cut}_2^{{\rm glob}}$ and $\gamma_2(t)$,  $0\leq t\leq T$,  be a shortest arc in $(SO_0(2,1)\times\mathbb{R},\rho_2)$ joining ${\rm Id}$ and $g$.
By virtue of Proposition 15, there exists only one shortest arc  
$\tilde{\gamma}_2(t)=\tilde{\gamma}_2(\varphi_0,\alpha_2,\beta;t)$, $0\leq t\leq T$, in $(SU(1,1)\times\mathbb{R},d_2)$ such that $\tilde{\gamma}_2(0)={\rm Id}$ and $\tilde{\Pi}(\tilde{\gamma}_2(t))=\gamma_2(t)$ for every $t\in [0,T]$. Put
$\tilde{\gamma}_2(T)=(A,B,v)$. Since ${\rm Re}(A)=0$ then $\beta\neq 0$. It follows from (\ref{form1}), (\ref{fform2}) that 
$$\tilde{\gamma}_2(\pi+\varphi_0+\beta T,\alpha_2,-\beta;T)=(\bar{A},-B,v)=(-A,-B,v).$$
It follows from here, (\ref{formula}), and (\ref{soot}) that
$\tilde{\Pi}(\tilde{\gamma}_2(\pi+\varphi_0+\beta T,\alpha_2,-\beta;T))=\gamma_2(T)=g,$
i.e., the points ${\rm Id}$ and $g$ in  $(SO_0(2,1)\times\mathbb{R},\rho_2)$ are connected by two different shortest paths parametrized by arclength. By Proposition 8, $g\in {\rm Cut}_2$.

It follows from Theorem 3 and Propositions 8, 16 that $(E,v)\notin {\rm Cut}_2$, where $v\in\mathbb{R}$, $E$ is the identity $3\times 3$--matrix. We assume that for some 
$g\notin {\rm Cut}_2^{{\rm loc}}\cup {\rm Cut}_2^{{\rm glob}}$ in $(SO_0(2,1)\times\mathbb{R},\rho_2)$
there exist two distinct shortest arcs
$\tilde{\Pi}(\tilde{\gamma}_2(\varphi_0,\alpha_2,\beta;t))$ and $\tilde{\Pi}(\tilde{\gamma}^{\ast}_2(\varphi^{\ast}_0,\alpha_2^{\ast},\beta^{\ast};t))$, $0\leq t\leq T$, joining  ${\rm Id}$ and $g$.
Let $\tilde{\gamma}_2(\varphi_0,\alpha_2,\beta;T):=(A,B,v)$ and $\tilde{\gamma}^{\ast}_2(\varphi^{\ast}_0,\alpha_2^{\ast},\beta^{\ast};T):=(A^{\ast},B^{\ast},v^{\ast})$.
By (\ref{formula}), we have $(A,B,v)=(A^{\ast},B^{\ast},v^{\ast})$ or $(A,B,v)=(-A^{\ast},-B^{\ast},v^{\ast})$. 
In the first case
$g\in {\rm Cut}_2^{{\rm loc}}$,  which is impossible. In the second case without loss of generality, we can assume that ${\rm Re}(A)<0$. Since  $\tilde{\gamma}_2(\varphi_0,\alpha_2,\beta;t)$, $0\leq t\leq T$,
connects ${\rm Id}$,  $(A,B,v)$ and ${\rm Re}(A)<0$, then there exists $t_1\in (0,T)$ such that
${\rm Re} (A_1)=0$, where $\tilde{\gamma}_2(\varphi_0,\alpha_2,\beta;t_1)=(A_1,B_1,v_1)$. Then $\tilde{\Pi}(\tilde{\gamma}_2(\varphi_0,\alpha_2,\beta;t_1))\in {\rm Cut}_2^{\rm glob}$.  This contradicts the fact that $\tilde{\Pi}(\tilde{\gamma}_2(\varphi_0,\alpha_2,\beta;t))$, with $0\leq t\leq T$, is a shortest arc.  
\noindent $\square$

\vspace{2mm}

The proof of the following proposition is similar to the proof of Proposition 11; we will omit it.

\vspace{2mm}

{\bf Proposition 18.} {\sl
	Let $\gamma_1(\varphi_0,\alpha_2,\beta;t)$, $0\leq t\leq T$, be a noncontinuable shortest arc in
	$(SO_0(2,1)\times\mathbb{R},\rho_1)$.
	Then $\gamma_2(\varphi_0,\alpha_2,\beta+\alpha_2;t)$, $0\leq t\leq T$, is a noncontinuable shortest arc in $(SO_0(2,1)\times\mathbb{R},\rho_2)$. The converse is also true.}

\vspace{2mm}

{\bf Proposition 19.} {\sl
	The cut locus ${\rm Cut}_1$ of the sub-Riemannian space $(SO_0(2,1)\times\mathbb{R},\rho_1)$ corresponding to ${\rm Id}$ is the set
	$${\rm Cut}_1={\rm Cut}_1^{{\rm loc}}\cup{\rm Cut}_1^{{\rm glob}},$$
	where 
	$${\rm Cut}_1^{{\rm loc}}=
	\left\{(C,v)\mid\,
	C=\left(
	\begin{array}{ccc}
	1 & 0 & 0 \\
	0 & \cos{\psi} & -\sin{\psi} \\
	0 & \sin{\psi} & \cos{\psi}
	\end{array}
	\right),\,\,v\in\mathbb{R},\,\,\psi\neq v+2\pi n,\,\,n\in\mathbb{Z}
	\right\};$$
	$${\rm Cut}_1^{{\rm glob}}=
	\left\{(C^{\ast},v)\mid\ C^{\ast}=\left(
	\begin{array}{ccc}
	c_{11} & c_{12}\cos{v}+c_{13}\sin{v} & c_{13}\cos{v}-c_{12}\sin{v} \\
	c_{21} & c_{22}\cos{v}+c_{23}\sin{v} & c_{23}\cos{v}-c_{22}\sin{v} \\
	c_{31} & c_{32}\cos{v}+c_{33}\sin{v} & c_{33}\cos{v}-c_{32}\sin{v}
	\end{array}\right),\right.$$
	$$\left.\quad C=(c_{ij})\in SO_0(2,1),\,\,c_{21}=-c_{12}, \,\, c_{31}=-c_{13},\,\,c_{23}=c_{32},\,\,c_{22}+c_{33}<-2,\,\,v\in\mathbb{R}
	\right\}.$$}

{\sl Proof.}
It follows from the definition of the cut set, Propositions 2, 18, and Theorem 4 that $g_2=(C,v)\in{\rm Cut}_2$ if and only if
$g_1=(C,v)\exp(-vE_3)=(C^{\ast},v)\in {\rm Cut}_1$, где
$$C^{\ast}=\left(
\begin{array}{ccc}
c_{11} & c_{12}\cos{v}+c_{13}\sin{v} & c_{13}\cos{v}-c_{12}\sin{v} \\
c_{21} & c_{22}\cos{v}+c_{23}\sin{v} & c_{23}\cos{v}-c_{22}\sin{v} \\
c_{31} & c_{32}\cos{v}+c_{33}\sin{v} & c_{33}\cos{v}-c_{32}\sin{v}
\end{array}
\right).$$
From here and Proposition 17 it follows that $g_2\in{\rm Cut}_2^{{\rm loc}}$ ($g_2\in{\rm Cut}_2^{{\rm glob}}$)  if and only if  $g_1\in {\rm Cut}_1^{{\rm loc}}$ (respectively $g_1\in {\rm Cut}_1^{{\rm glob}}$).
\noindent $\square$

\section{Noncontinuable shortest arcs on the Lie group $SO_0(2,1)\times\mathbb{R}$}

By Proposition 4, to find shortest arcs in $(SO_0(2,1)\times\mathbb{R},\rho_i)$, $i=1,2$, it suffices study the geodesic segments $\gamma_i(t)$, $0\leq t\leq T$, such that $\gamma_i(0)={\rm Id}$.

\vspace{2mm}

{\bf Theorem 10.} {\sl
	Let $\alpha_2\neq\pm 1$, $\beta\neq 0$ and  $\gamma_2(t)=\gamma_2(\varphi_0,\alpha_2,\beta;t)$, $0\leq t\leq T$, be a noncontinuable shortest arc of the sub-Riemannian space  $(SO_0(2,1)\times\mathbb{R},\rho_2)$. Then
	
	1. If $0<\frac{|\beta|}{\sqrt{1-\alpha_2^2}}<1$ then $T\in\left(\frac{\pi}{|\beta|},\frac{2\pi}{|\beta|}\right)$ and $T$ 
	satisfies the system of equations
	$$\cos\frac{|\beta|T}{2} =-\frac{|\beta|\sinh(T\sqrt{1-\alpha_2^2-\beta^2}/2)}{\sqrt{(1-\alpha_2^2)\cosh^2(T\sqrt{1-\alpha_2^2-\beta^2}/2)-\beta^2}},$$
	$$\sin\frac{|\beta|T}{2} =\frac{\sqrt{1-\alpha_2^2-\beta^2}\cosh(T\sqrt{1-\alpha_2^2-\beta^2}/2)}{\sqrt{(1-\alpha_2^2)\cosh^2(T\sqrt{1-\alpha_2^2-\beta^2}/2)-\beta^2}}.$$
	
	2. If $\frac{|\beta|}{\sqrt{1-\alpha_2^2}}=1$ then $T\in\left(\frac{\pi}{|\beta|},\frac{2\pi}{|\beta|}\right)$ and $T$ satisfies the system of equations
	$$\cos\frac{\beta T}{2}=-\frac{\beta T/2}{\sqrt{1+\beta^2T^2/4}},\quad \sin\frac{\beta T}{2}=\frac{1}{\sqrt{1+\beta^2T^2/4}}.$$
	
	3. If $1<\frac{|\beta|}{\sqrt{1-\alpha_2^2}}<\frac{2}{\sqrt{3}}$ then  $T\in\left(\frac{\pi}{|\beta|},\frac{2\pi}{|\beta|}\right)$ and $T$ satisfies the system of equations
	\begin{equation}
	\label{soot1}
	\cos\frac{|\beta|T}{2} =-\frac{|\beta|\sin(T\sqrt{\alpha_2^2+\beta^2-1}/2)}{\sqrt{\beta^2-(1-\alpha_2^2)\cos^2(T\sqrt{\alpha_2^2+\beta^2-1}/2)}},
	\end{equation}
	\begin{equation}
	\label{soot2}
	\sin\frac{|\beta|T}{2} =\frac{\sqrt{\alpha_2^2+\beta^2-1}\cos(T\sqrt{\alpha_2^2+\beta^2-1}/2)}
	{\sqrt{\beta^2-(1-\alpha_2^2)\cos^2(T\sqrt{\alpha_2^2+\beta^2-1}/2)}}.
	\end{equation}
	
	4. If $\frac{|\beta|}{\sqrt{1-\alpha_2^2}}=\frac{2}{\sqrt{3}}$ then $T=\frac{2\pi}{|\beta|}$.
	
	5. If $\frac{2}{\sqrt{3}}<\frac{|\beta|}{\sqrt{1-\alpha_2^2}}<\frac{3}{\sqrt{5}}$ then
	$T\in\left(\frac{2\pi}{|\beta|},\frac{3\pi}{|\beta|}\right)$ and $T$ satisfies the system of equations
	(\ref{soot1}), (\ref{soot2}).
	
	6. If $\frac{|\beta|}{\sqrt{1-\alpha_2^2}}\geq\frac{3}{\sqrt{5}}$ then $T=\frac{2\pi}{\sqrt{\beta^2+\alpha_2^2-1}}$.
}

\vspace{2mm}

{\sl Proof.}
Let $\alpha_2\neq\pm 1$, $\beta\neq 0$ and $\gamma_2(t)=\gamma_2(\varphi_0,\alpha_2,\beta;t)$,  $0\leq t\leq T$, be a noncontinuable
shortest arc. By Proposition 17,   
$g=\gamma_2(T)$   belongs to the union of ${\rm Cut}_2^{{\rm loc}}$ and ${\rm Cut}_2^{{\rm glob}}$. 
Then, by Theorem 3, we get

1) if $g\in {\rm Cut}_2^{{\rm loc}}$ then $\frac{|\beta|}{\sqrt{1-\alpha_2^2}}>1$ and $T=\frac{2\pi}{\sqrt{\beta^2+\alpha_2^2-1}}$ due to (\ref{mn01})--(\ref{mn03});

2) if  $g\in {\rm Cut}_2^{{\rm glob}}$ then
\begin{equation}
\label{lji}
n_2(T)\cos\frac{\beta T}{2}+\beta m_2(T)\sin\frac{\beta T}{2}=0.
\end{equation}

Let $T_0$ be the smallest positive root of equality (\ref{lji}). 
If $0<\frac{|\beta|}{\sqrt{1-\alpha_2^2}}\leq 1$ then $T=T_0$. If $\frac{|\beta|}{\sqrt{1-\alpha_2^2}}>1$ then $T={\rm min}\left\{T_0,\frac{2\pi}{\sqrt{\beta^2+\alpha_2^2-1}}\right\}$.

Let $0<\frac{|\beta|}{\sqrt{1-\alpha_2^2}}<1$. 
Define the function $F_1(t)$ on the segment $\left[0,\frac{2\pi}{|\beta|}\right]$ by the
formula
$$F_1(t)=\cosh\frac{t\sqrt{1-\alpha_2^2-\beta^2}}{2}\cos\frac{{|\beta|t}}{2}+\frac{|\beta|}{\sqrt{1-\alpha_2^2-\beta^2}}\sinh\frac{t\sqrt{1-\alpha_2^2-\beta^2}}{2}\sin\frac{{|\beta|t}}{2}.$$
By virtue of (\ref{mn02}), $F_1(t)$ is the left side of the equation (\ref{lji}). 
It is easy to see that
$$F_1(0)=1,\quad F_1\left(\frac{2\pi}{|\beta|}\right)<0,\quad F_2^{\prime}(t)=\frac{(1-\alpha_2^2)\cos\frac{|\beta|t}{2}\sinh\frac{t\sqrt{1-\alpha_2^2-\beta^2}}{2}}{2\sqrt{1-\alpha_2^2-\beta^2}}.$$
The function  $F_1(t)$ increases on $\left[0,\frac{\pi}{|\beta|}\right]$, decreases on $\left[\frac{\pi}{|\beta|},\frac{2\pi}{|\beta|}\right]$, and has the unique root  $T_0\in \left(\frac{\pi}{|\beta|},\frac{2\pi}{|\beta|}\right)$. 
This implies item~1 of Theorem 10.

Let  $\frac{|\beta|}{\sqrt{1-\alpha_2^2}}=1$. By virtue of (\ref{mn01}),  the equation (\ref{lji})
will be written as
$$\cos\frac{|\beta|T}{2}+\frac{|\beta| T}{2}\sin\frac{|\beta|T}{2}=0.$$
It is easy to see that the smallest positive root of this equation belongs to the interval
$\left(\frac{\pi}{|\beta|},\frac{2\pi}{|\beta|}\right)$. 
This implies item~2 of Theorem 10.

Let  $\frac{|\beta|}{\sqrt{1-\alpha_2^2}}>1$. 
Define the function $F_2(t)$ on the segment $\left[0,\frac{2\pi}{\sqrt{\beta^2+\alpha_2^2-1}}\right]$
by the formula
\begin{equation}
\label{ll}
F_2(t)=\cos\frac{t\sqrt{\beta^2+\alpha_2^2-1}}{2}\cos\frac{{|\beta|t}}{2}+\frac{|\beta|}{\sqrt{\beta^2+\alpha_2^2-1}}\sin\frac{t\sqrt{\beta^2+\alpha_2^2-1}}{2}\sin\frac{{|\beta|t}}{2}.
\end{equation}
By virtue of (\ref{mn03}), $F_2(t)$ is the left side of the equation (\ref{lji}). 
It is easy to see that $F_2(0)=1>0$,
$$F_2\left(\frac{2\pi}{|\beta|}\right)=-\cos\frac{\pi\sqrt{\alpha_2^2+\beta^2-1}}{|\beta|} ,\quad  F_2^{\prime}(t)=\frac{(1-\alpha_2^2)\cos\frac{|\beta|t}{2}\sin\frac{t\sqrt{\alpha_2^2+\beta^2-1}}{2}}{2\sqrt{\alpha_2^2+\beta^2-1}}.$$  
The function  $F_2(t)$ increases on the segments $\left[0,\frac{\pi}{|\beta|}\right]$, 
$\left[\frac{2\pi}{|\beta|},{\rm min}\left\{\frac{3\pi}{|\beta|},\frac{2\pi}{\sqrt{\alpha_2^2+\beta^2-1}}\right\}\right]$ and decreases on $\left[\frac{\pi}{|\beta|},\frac{2\pi}{|\beta|}\right]$.

If $0<\frac{|\beta|}{\sqrt{1-\alpha_2^2}}<\frac{2}{\sqrt{3}}$ then 
$F_2\left(\frac{2\pi}{|\beta|}\right)<0$. Therefore, the smallest positive root $T_0$ of the equation (\ref{lji}) belongs to the interval $\left(\frac{\pi}{|\beta|},\frac{2\pi}{|\beta|}\right)$ and 
$T=T_0$. This and (\ref{ll}) imply item~3 of Theorem 10.

If $\frac{|\beta|}{\sqrt{1-\alpha_2^2}}=\frac{2}{\sqrt{3}}$ then
$F_2\left(\frac{2\pi}{|\beta|}\right)=0$. This implies item~4 of Theorem 10.

If $\frac{|\beta|}{\sqrt{1-\alpha_2^2}}>\frac{2}{\sqrt{3}}$ then
$F_2\left(\frac{2\pi}{|\beta|}\right)>0$, and the function $F_2(t)$ is positive on the segment $\left[0,\frac{2\pi}{|\beta|}\right]$.  If at the same time
$\frac{2\pi}{\sqrt{\alpha_2^2+\beta^2-1}}\leq\frac{3\pi}{|\beta|}$, i.e., 
$\frac{|\beta|}{\sqrt{1-\alpha_2^2}}\geq \frac{3}{\sqrt{5}}$, then $F_2(t)>0$ on the segment $\left[0,\frac{2\pi}{\sqrt{\beta^2+\alpha_2^2-1}}\right]$. Thus $T=\frac{2\pi}{\sqrt{\alpha_2^2+\beta^2-1}}$, and item~6 of Theorem 10 is proved.

If $\frac{2\pi}{\sqrt{\alpha_2^2+\beta^2-1}}>\frac{3\pi}{|\beta|}$, i.e.,
$\frac{2}{\sqrt{3}}<\frac{|\beta|}{\sqrt{1-\alpha_2^2}}<\frac{3}{\sqrt{5}}$, then 
$F_2\left(\frac{3\pi}{|\beta|}\right)=-\cos\frac{\pi |\beta|}{\sqrt{\alpha_2^2+\beta^2-1}}<0$, and $T_0\in\left(\frac{2\pi}{|\beta|},\frac{3\pi}{|\beta|}\right)$. This and (\ref{ll}) imply item~5 of Theorem 10.
\noindent $\square$

\vspace{2mm}

By virtue of Proposition 18, if $\gamma_1(\varphi_0,\alpha,\beta,t)$, $0\leq t\leq T$, is a noncontinuable shortest arc in $(SO_0(2,1)\times\mathbb{R},\rho_1)$, then $T$ satisfies items 1--6 of Theorem 7 after replacing $\beta$ by $\beta+\alpha_2$. 

The following is immediate from Propositions 2, 18.

\vspace{2mm}

{\bf Proposition 20.} {\sl
	For every $(C,v)\in SO_0(2,1)\times\mathbb{R}$,
	$\rho_1({\rm Id},(C,v))=\rho_2({\rm Id},(\tilde{C},v)),$
	where
	$$\tilde{C}=\left(
	\begin{array}{cccc}
	c_{11} & c_{12}\cos{v}-c_{13}\sin{v} & c_{13}\cos{v}+c_{12}\sin{v} \\
	c_{21} & c_{22}\cos{v}-c_{23}\sin{v} & c_{23}\cos{v}+c_{22}\sin{v} \\
	c_{31} & c_{32}\cos{v}-c_{33}\sin{v} & c_{33}\cos{v}+c_{32}\sin{v}
	\end{array}
	\right)$$ }

The following is immediate from Theorem 7 and (\ref{abc2}).

\vspace{2mm}

{\bf Proposition 21.} {\sl For every $(C,v)\in SO_0(2,1)\times\mathbb{R}$,
	$\rho_2^2({\rm Id},(C,v))=v^2+\rho_2^2({\rm Id},(C,0))$.}

\vspace{2mm}

Exact formulas for $\rho_2({\rm Id},(C,0))$ are found in [7]. Using them and Propositions 20, 21, 
we can obtain exact formulas for $\rho_i({\rm Id},(C,v))$, $i=1,2$. We will not
present them because of their bulkiness.

\vspace{2mm}

{\bf Acknowledgments.}
The author thanks Professor V.N. Berestovskii for useful discussions and remarks.

\end{document}